\PassOptionsToPackage{unicode}{hyperref}
\PassOptionsToPackage{hyphens}{url}
\documentclass[article]{article}
\usepackage{amsmath,amssymb}
\usepackage{iftex}
\ifPDFTeX
  \usepackage[T1]{fontenc}
  \usepackage[utf8]{inputenc}
  \usepackage{textcomp} 
\else 
  \usepackage{unicode-math} 
  \defaultfontfeatures{Scale=MatchLowercase}
  \defaultfontfeatures[\rmfamily]{Ligatures=TeX,Scale=1}
\fi
\usepackage{lmodern}
\ifPDFTeX\else
\fi
\IfFileExists{upquote.sty}{\usepackage{upquote}}{}
\IfFileExists{microtype.sty}{
  \usepackage[]{microtype}
  \UseMicrotypeSet[protrusion]{basicmath} 
}{}
\makeatletter
\@ifundefined{KOMAClassName}{
  \IfFileExists{parskip.sty}{%
    \usepackage{parskip}
  }{
    \setlength{\parindent}{0pt}
    \setlength{\parskip}{6pt plus 2pt minus 1pt}}
}{
  \KOMAoptions{parskip=half}}
\makeatother
\usepackage{xcolor}
\usepackage[margin=1in]{geometry}
\usepackage{graphicx}
\makeatletter
\def\maxwidth{\ifdim\Gin@nat@width>\linewidth\linewidth\else\Gin@nat@width\fi}
\def\maxheight{\ifdim\Gin@nat@height>\textheight\textheight\else\Gin@nat@height\fi}
\makeatother
\setkeys{Gin}{width=\maxwidth,height=\maxheight,keepaspectratio}
\makeatletter
\def\fps@figure{htbp}
\makeatother
\providecommand{\tightlist}{%
  \setlength{\itemsep}{0pt}\setlength{\parskip}{0pt}}
\usepackage{graphics}
\usepackage{amsmath, amsfonts, amsthm, amssymb, amscd,a4wide}
\usepackage{algorithmic}
\usepackage[linesnumbered,ruled,vlined]{algorithm2e}
\usepackage{hyperref}
\usepackage{url}
\usepackage{float}
\usepackage{epstopdf}
\usepackage{subcaption}
\let\oldmarginpar\marginpar
\renewcommand\marginpar[1]{\-\oldmarginpar[\raggedleft\footnotesize #1]%
{\raggedright\footnotesize #1}}

\floatname{algorithm}{Algorithm}

\usepackage{tikz} 
\usepackage{xcolor}
\hypersetup{frenchlinks=true}

\numberwithin{equation}{section}

\newcommand \be   {\begin{equation}}
\newcommand \bel {\begin{equation}\label}
\newcommand \ee   {\end{equation}}

\newcommand \RR    {\mathbb{R}}

\newcommand \la         \langle
\newcommand \ra     \rangle

\usepackage{mathrsfs}

\usepackage{authblk}

\tikzstyle{startstop} = [rectangle, rounded corners, minimum width=3cm, minimum height=1cm,text centered, draw=black ]
\tikzstyle{io} = [rectangle, rounded corners, minimum width=3cm, minimum height=1cm,text centered, draw=black] 
\tikzstyle{process} = [rectangle, minimum width=3cm, minimum height=1cm, text centered, draw=black, fill=orange!30]
\tikzstyle{decision} = [diamond, minimum width=3cm, minimum height=1cm, text centered, text - white, draw=black, fill=black!30] 
\tikzstyle{arrow} = [thick,->,>=stealth]
\tikzstyle{bbox} = [rectangle, minimum width=3cm, minimum height=1cm, text centered, text = white, draw=black, fill=black]
\tikzstyle{pp} = [rectangle, minimum width=3cm, minimum height=0.5cm, text centered, draw=black] 
\tikzstyle{pp1} = [rectangle, draw=black!50, thick, minimum width=0.5cm, minimum height = 0.5cm]
\tikzstyle{crl} = [circle, draw=black!50, thick, minimum size = 1.5cm]
\tikzstyle{crl1} = [circle, draw=black!50, thick, minimum size = 0.7cm]
\tikzstyle{line} = [draw, -latex']
\newsavebox{\tempbox}

\tikzstyle{block} = [draw, rectangle, 
    minimum height=3em, minimum width=6em]
\tikzstyle{sum} = [draw, fill=blue!20, circle, node distance=1cm]
\tikzstyle{input} = [coordinate]
\tikzstyle{output} = [coordinate]
\tikzstyle{pinstyle} = [pin edge={to-,thin,black}]
\usepackage{booktabs}
\usepackage{longtable}
\usepackage{array}
\usepackage{multirow}
\usepackage{wrapfig}
\usepackage{float}
\usepackage{colortbl}
\usepackage{pdflscape}
\usepackage{tabu}
\usepackage{threeparttable}
\usepackage{threeparttablex}
\usepackage[normalem]{ulem}
\usepackage{makecell}
\usepackage{xcolor}
\ifLuaTeX
  \usepackage{selnolig}  
\fi
\usepackage{bookmark}
\IfFileExists{xurl.sty}{\usepackage{xurl}}{} 
\hypersetup{
  pdftitle={A class of kernel-based scalable algorithms for data science},
  pdfauthor={Philippe G. LeFloch, Jean-Marc Mercier\^{}\textbackslash dag, and Shohruh Miryusupov},
  hidelinks,
  pdfcreator={LaTeX via pandoc}}

\title{A class of kernel-based scalable algorithms for data science}
\author{Philippe G.
LeFloch\footnote{Laboratoire Jacques-Louis Lions, Sorbonne University and Centre National de la Recherche Scientifique, 4 Place Jussieu, 75258 Paris, France. Email: contact@philippelefloch.org},
Jean-Marc Mercier\(^\dag\), and Shohruh
Miryusupov\footnote{MPG-Partners, 166 rue de Stockholm, 75008 Paris, France. Email: jean-marc.mercier@mpg-partners.com, :shohruh.miryusupov@mpg-partners.com \ First version: October 2024. Revised version: December 2024.}}
\date{}

\begin{document}
\maketitle
\begin{abstract}
We present several generative and predictive algorithms based on the
RKHS (reproducing kernel Hilbert spaces) methodology, which, most
importantly, are scale up efficiently with large datasets or
high-dimensional data. It is well recognized that the RKHS methodology
leads one to efficient and robust algorithms for numerous tasks in data
science, statistics, and scientific computation. However, the
implementations existing the literature are often difficult to scale up
for encompassing large datasets. In this paper, we introduce a simple
and robust, divide-and-conquer methodology. It applies to large scale
datasets and relies on several kernel-based algorithms, which
distinguish between various extrapolation, interpolation, and optimal
transport steps. We argue how to select the suitable algorithm in
specific applications thanks to a feedback of perfomance criteria. Our
primary focus is on applications and problems arising in industrial
contexts, such as generating meshes for efficient numerical simulations,
designing generators for conditional distributions, constructing
transition probability matrices for statistical or stochastic
applications, and addressing various tasks relevant to the Artificial
Intelligence community. The proposed algorithms are highly relevant to
supervised and unsupervised learning, generative methods, as well as
reinforcement learning.
\end{abstract}

\setcounter{secnumdepth}{2}
\setcounter{tocdepth}{2} 
\tableofcontents

\section{Introduction}\label{introduction}

Reproducing Kernel Hilbert Space (RKHS) methods provide a robust and
flexible framework for regression, positioning them as a compelling
alternative to widely used techniques such as neural networks and
tree-based models. These methods excel in both interpolation and
extrapolation tasks and play a pivotal role in the development of
optimal transportation algorithms, making them applicable across diverse
domains. In Artificial Intelligence, RKHS-based approaches significantly
contribute to supervised and unsupervised learning, generative models,
reinforcement learning, and numerous areas of applied mathematics,
including numerical simulations, partial differential equations, and
stochastic analysis. These versatile capabilities underscore the
importance of RKHS-based approaches in modern data science, financial
mathematics, and statistics.

Kernel methods exhibit several distinct advantages over other
approaches. The first key advantage is their reproducibility, which
enables them to exactly replicate arbitrary discrete function values---a
feature particularly valuable in specific applications. Moreover, kernel
methods are regarded as white-box models, providing transparency by
offering worst-case error bounds for both extrapolation and
interpolation tasks. These error bounds are crucial for performing error
analysis in numerical simulations and for calculating confidence
intervals in statistical applications. Furthermore, as this paper will
demonstrate, worst-case error bounds serve as powerful tools for
optimizing kernel methods, thereby enhancing their efficiency and
reliability.

Despite their many advantages, kernel methods face significant
challenges when applied to large-scale datasets due to their inherent
algorithmic complexity, which often renders them impractical. After
briefly reviewing kernel-based methods, we introduce (in Section
\ref{multi-scale-kernel-methods}) a comprehensive strategy for handling
large-scale datasets, referred to as the multi-scale kernel method. This
approach aligns with strategies proposed in active research areas, such
as the Nyström method, low-rank approximations, and clustering
techniques, and shares similarities with works like \cite{Falkon:2020}
and \cite{BLWY2}. The proposed multi-scale kernel method retains key
advantages of kernel-based approaches, such as reproducibility and error
estimation, even for large datasets, while significantly simplifying
implementation on distributed architectures. It is based on
computational units that leverage multi-CPU or GPU architectures, with
each unit operating within an independent memory space defined by a
clustering method. By employing robust divide-and-conquer techniques,
our strategy effectively addresses the challenges of large-scale
datasets, enabling kernel algorithms to adapt efficiently to diverse
computational environments with varying memory and computational
constraints.

This approach relies on the quality of an initial coarse regressor,
requiring scalable algorithms capable of performing the initial
partitioning at a lower computational cost than the overall prediction
process. While algorithms like k-means are fast, scalable, and natural
candidates for this task, they come with certain limitations: k-means is
optimized for a specific kernel (Euclidean distance) and typically
produces linear separations, making it sensitive to scale. Moreover, it
is prone to becoming trapped in local minima. To address these
shortcomings, this paper explores a range of alternative scalable
algorithms tailored to optimize kernel methods for large-scale datasets.
Our approach incorporates simple, yet original, and robust techniques,
selected for their demonstrated numerical efficiency in practice. These
methods are presented in Section \ref{general-purpose-algorithms} and
applied to the multi-scale approach in Section
\ref{clustering-algorithms}.

The same \textit{multi-scale} approach can be effectively applied to
large-scale optimal transport problems. In many applications, such as
generative models, it is crucial to define a mapping \(T\) between two
distributions \(X\) and \(Y\), typically represented in transport theory
as \(T_{\#}X = Y\), denoting the \textit{push-forward} of one
distribution to another. The multi-scale method is well-suited to this
context, as the problem can be divided across computational units for
efficient processing. Our experiments demonstrate that modern processor
cores can handle dataset sizes ranging from \(1,000\) to \(10,000\)
elements per second. As a result, in Section
\ref{A-framework-for-dealing-with-two-distributions}, we have tailored
our tools for small to medium-sized datasets and introduced key
algorithms based on this optimal transport framework.

Transportation problems have a long history and are closely tied to
operational research, with classical formulations such as the linear sum
assignment problem (LSAP), weighted graph matching problem (WGMP), the
traveling salesman problem (TSP), and graph partitioning (GP). These
problems can be traced back to the machine learning community as early
as the 1980s (see \cite{HT:1985}). Beginning with \cite{Cuturi:2013}, a
substantial body of literature has emerged (see \cite{SFVTP:2023} and
references therein), focusing on transportation problems within the
framework of optimal transport theory, particularly with entropy-based
regularization. A comprehensive review of the theory of optimal
transport can be found in \cite{Villani:2009}. The practical
implementations primarily employ the Sinkhorn-Knopp algorithm, a
standard descent method applied to Wasserstein or Gromov-Wasserstein
problems \cite{MF2011,Cuturi:2016,Bach} via the Kantorovich formulation,
which is the dual to Monge's original formulation.

Building on our earlier investigations in
\cite{PLF-JMM-estimate}-\cite{PLF-JMM-Wilmott}-\cite{PLF-JMM:2023} and
\cite{LMM-SSRN1}-\cite{LMM}, we adopt a slightly different approach here
and address distinct constraints and objectives. Specifically, we often
proceed without regularization, aiming to respect an exact transport
property by defining mappings that precisely satisfy \(T_{\#}X = Y\),
leveraging the reproducibility properties of kernel methods. Our focus
is on methods that, in principle, achieve the best possible convergence
rates, reducing the computational burden of kernel methods and
facilitating the successful passage of statistical tests. In this
context, we explore \textit{combinatorial algorithms} rather than
descent-based methods. Numerically, this is implemented through
\textit{pairing} or \textit{matching} algorithms, such as the linear sum
assignment problem (LSAP) or the traveling salesman problem (TSP),
instead of continuous descent algorithms like Sinkhorn-Knopp. This
direct, combinatorial approach is particularly effective for low- to
medium-sized datasets, our primary focus in the present work, and can be
parallelized to handle both regularized and non-regularized problems.
Finally, the two algorithmic paradigms---direct vs.~regularized, Monge
vs.~Kantorovich---are complementary: regularized problems typically rely
on descent algorithms, for which a combinatorial method can provide a
high-quality initial solution, thereby accelerating the descent process.

Our research led to the development of CodPy, an open-source
kernel-based library introduced in our companion textbook \cite{LMM}.
All numerical illustrations of the algorithms presented in this paper
are based on this library. Its open-source nature ensures that our
numerical experiments are fully reproducible and, in the future, it may
serve as a reliable benchmark for comparing the proposed algorithms with
alternative methods. For a detailed framework on RKHS methods, including
their application to designing \emph{physics-informed} methods and
solving partial differential equations, we refer the interested reader
to \cite{LMM,LMMprep}, as well as \cite{Bachdeux}, and the references
cited therein.

\section{Fundamentals on kernel
methods}\label{fundamentals-on-kernel-methods}

\subsection{Aim}\label{aim}

We begin with basic material on the mathematical foundations that are
required for the definition and study of the algorithms presented in
this paper. Reproducing kernel methods can be thought as universal
regressors, which enjoy optimal error bounds for dealing with standard
operations, such as the extrapolation of datasets. This main feature of
the RKHS methods will now be briefly explained. For a detailed
introduction to the theory of reproducing kernel Hilbert space (RKHS),
we refer the reader to \cite{BTA}. Additionally, reproducing kernel
methods can be applied to compute mappings that transport a distribution
into another, and this important issue will be discussed later on in
this text. Below, we also describe a multi-scale methodology and perform
extrapolation with such kernel methods.

\subsection{Framework for this paper}\label{framework-for-this-paper}

Kernels are employed to represent functions of the general form
\(x \mapsto f(x)+\epsilon(x)\), where \(\epsilon(x)\) an (optional)
white noise term. It is convenient to introduce the following matrix
notation for the variables and the vector-valued maps under
consideration: \begin{equation}
\label{X}
X :=\{x^n_d\}_{n,d=1}^{N_X,D} \in \RR^{N_X,D}, \quad f(X) :=\{f(x^n_d) \}_{n,d=1}^{N_X,D_f} \in \RR^{N_X,D_f}. 
\end{equation} Here, \(X\) denotes a set of points in \(\RR^D\), which
serves as the input to a continuous function \(f(\cdot) \in \RR^{D_f}\).
Observe that the set \((X,f(X))\) can also be interpreted as a discrete
probability distribution \(\frac{1}{N_X}\sum_n \delta_{x^n}\), in which
\(\delta_{x}\) is the Dirac mass. Interestingly, most of the algorithms
in the present paper could be extended to arbitrary (continuous)
probability distributions. In the language of supervised learning, the
set \((X,f(X))\) is referred to as the \textit{training set}.

By definiion, a function \(k: \RR^D \times \RR^D \mapsto \mathbb{R}\) is
called a \textit{kernel} if it is symmetric and positive definite (see
\cite{BTA} for a definition). Given two sets of points
\(X \in \RR^{N_X,D}\) and \(Y \in \RR^{N_Y,D}\), we introduce their
associated \textit{Gram matrix} \begin{equation}
K(X,Y):= \big(k(x^i,y^j)\big)_{i,j} \in \mathbb{R}^{N_X,N_Y}.
\end{equation} Positive-definite kernels corresopnd to symmetric,
positive-definite Gram matrices \(K(X,X)\), which are invertible
provided \(X\) consists of \textit{distinct} points in \(\RR^{N_X,D}\).
Given a set \(X\), a given kernel can be used to construct a space of
vector-valued functions (referred to as regressors) parameterized by
matrices in \(\mathbb{R}^{N_Y,D_\theta}\), as follows: \begin{equation}
\label{Hk}
\mathcal{H}_k^Y := \Big\{f_{k,\theta}(\cdot) := \sum_{n=1}^{N_Y} \theta^n k(\cdot,y^n) = K(\cdot,Y) \theta, \quad \theta \in \mathbb{R}^{N_Y,D_\theta} \Big\}.
\end{equation} Here, the parameters denoted by \(\theta\) are computed,
or \textit{fitted}, to a given continuous function, defined on the
training set \(X,f(X)\), according to the relation
\begin{equation}\label{FIT}
 f_{k, \theta} (\cdot) := K(\cdot,Y) \theta, \quad \theta = \Big( K(X,Y) + \epsilon R(X,Y) \Big)^{-1} f(X),
\end{equation} in which \(\epsilon \ge 0\) and \(R(X,Y)\) determine an
optional regularizing term.

It is important to observe that the right-hand side contains a
\textit{rectangular} matrix if \(N_Y < N_X\), but a \textit{square}
matrix if \(N_X=N_Y\) ---in particular if \(Y = X\) coincides.
Typically, the implied system is solved using a least-square approach,
although other methods may be applied depending on the specific problem
at hand. The overall number of operations, that is, the number of
elementary operations to perform a computation \eqref{FIT} is
\(\mathcal{O}(N_X^2 N_Y+N_Y^3)\), thus the primary role of the set \(Y\)
is to reduce the computational complexity of the matrix inversion. The
algorithms associated with the computation of \(Y\) are related to a
notion of \textit{clustering}, which will be presented in Section
\ref{clustering-algorithms}, below.

The fitting procedure defined by the equation \eqref{FIT} can be
expressed in terms of the following operator \(\mathcal{P}_k\):
\begin{equation}
 f_{k, \theta}(\cdot) = \mathcal{P}_k f(X), \qquad \mathcal{P}_k =K(\cdot,Y)\big( K(X,Y) + \epsilon R(X,Y) \big)^{-1}.
\end{equation} This operator standpoint has several advantages, as it
conveniently allows for the (matrix) definition and for handling many
other common operators. For example, the kernel-based
\textit{gradient operator} can be defined by \begin{equation}
\label{nabla}
 \nabla f_k(\cdot,\theta) = \nabla_k f(X), \quad \nabla_k = (\nabla K)(\cdot,Y)\big( K(X,Y) + \epsilon R(X,Y) \big)^{-1}
\end{equation} and, similarly, we can design a discrete version of the
Laplace-Beltrami operator by \begin{equation}
\label{delta}
  \Delta_k = -\nabla_k^T \nabla_k.
\end{equation}

We refer to the notions with \(Y=X\) and without a regularization
(\(\epsilon=0\)) as the \textit{extrapolation procedure}, while the
other choices are referred to as \textit{interpolation procedures.} For
simplicitly, let us first focus on this purely extrapolation procedure,
which is also called here the extrapolation mode and exhibits many
interesting features and, however, showcases an overall algorithmic
burden of \(\mathcal{O}(N_X^3)\) operations ---a prohibitive cost to
handle large scale datasets. In this context, kernel regressors provide
us with error estimates of the form \begin{equation}\label{dkerr}
 \| f(\cdot)-f_{k, \theta}(\cdot) \|_{\ell^2} \le d_k(\cdot,X) \,  \|f\|_{\mathcal{H}_k^{\mathbb{X}}}. 
\end{equation} We refer to \cite{BTA} for the precise notions. Roughly
speaking, \(\mathcal{H}_k^{\mathbb{X}}\) denotes the functional space
associated with the discrete space \eqref{Hk} generated by the kernel
\(k\), localized by the measure \(\mathbb{X}\) to which the set of
points \(X\) is supposed to be i.i.d. (independent and identically
distributed) random variables. The factor
\(\|f\|_{\mathcal{H}_k^{\mathbb{X}}}\) can be approximated by the lower
bound \(<f(X),\theta>_{\ell^2}\). The factor \(d_k(\cdot,X)\) is the
so-called kernel discrepancy distance. Numerically, this distance is
evaluated on a discrete set of points \(Z \in \RR^{N_Z,D}\), which in
the context of machine learning can be referred to as the
\textit{test set}. A discrete version of this distance, also known as
the Maximum Mean Discrepancy (MMD), was introduced in \cite{GBRS} and
reads \begin{equation}
\label{dk}
    d_k\big(Z,X\big)^2 := \frac{1}{N_X^2} \sum_{n,m=1}^{N_x}  k(x^n,x^m) + \frac{1}{N_Z^2} \sum_{n,m}^{N_Z} k(z^n,z^m) - \frac{2}{N_X N_Z} \sum_{n,m=1}^{N_X,N_Z} k(x^n,z^m). 
\end{equation} This distance can also be generalized and defines a
distance \(d_k(\mathbb{Y},\mathbb{X})\) between any two probability
distributions \(\mathbb{X},\mathbb{Y}\). For any positive-definite
kernel, this distance is positive, and we also observe that it can be
generated from the formula \(d_k(x,z)=k(x,x)+k(z,z)-2k(x,z)\).

Let us also introduce the \textit{discrepancy matrix} induced by the
kernel \(k\), namely \begin{equation}\label{Mk}
 M_k(X,Z) := \big(d_k(x^n,z^m) \big)_{n,m=1}^{N_X,N_Z}. 
\end{equation} An interesting property of kernel extrapolation operators
is their ability to achieve zero-error on the training set in the
extrapolation mode, a property referred to as \textit{reproductibility},
that is, \begin{equation}\label{REPRO}
 f(X)=f_{k, \theta}(X),\quad \theta = K(X,X)^{-1}f(X). 
\end{equation} This property follows from \(\mathcal{P}_k(X,X) = I_d\),
or \eqref{dkerr}, by observing that \(d_k\big(X,X\big) = 0\).

\subsection{Dealing with two
distributions}\label{dealing-with-two-distributions}

\label{distrib}

In many applications, it is necessary to define mappings between two
distributions, say \(\mathbb{X}\) and \(\mathbb{Y}\), with support in
\(\mathbb{R}^{D_X}\) and \(\mathbb{R}^{D_Y}\), respectively. These
distributions are typically only known through two i.i.d. samples of
equal length, say \(X \in \RR^{N,D_X}\) and \(Y \in \RR^{N,D_Y}\). In
order to define such mappings, our construction is essentially based on
the theory of optimal transport.

We begin by recalling some key concepts from the theory of optimal
transport. A map \(T : \mathbb{R}^{D_X} \mapsto \mathbb{R}^{D_Y}\) is
said to \textit{transport} the distribution \(\mathbb{X}\) to the
distribution \(\mathbb{Y}\) if \(\mathbb{Y} = T \circ \mathbb{X}\) (in
the sense of composition of measures). In optimal transport theory, such
mappings are referred to as \textit{push forward maps} and are denoted
by \(T_{\#} \mathbb{X} = \mathbb{Y}\). From a discrete point of view,
and in thanks to \eqref{FIT} in the pure extrapolation mode (with
\(\epsilon=0\)), the map \[
  Y_{k, \theta}(\cdot) := K(\cdot,X) \theta, \quad \theta = K(X,X)^{-1} Y
\] defines a push forward of the discrete distribution \(\delta_X\) to
\(\delta_Y\), where \(\delta_X := \frac{1}{N}\sum_n \delta_{x^n}\),
\(\delta_{x}\) is the Dirac measure centered at \(x\). It is important
to observe, however, that all index permutations
\(\sigma: [1,...,N] \mapsto [1,...,N]\) of the set \(Y\) also define
push-forward maps. We are usually interested in finding the permutations
that define smooth, invertible transport maps, which take the form
\begin{equation} 
\label{ED}
  (Y\circ \sigma)_{k,\theta}(\cdot) := K(\cdot,X) \theta, \quad \theta = K(X,X)^{-1} (Y\circ \sigma). 
\end{equation} Hence, from the discrete standpoint, finding a smooth and
invertible transport map from \(X\) to \(Y\) is equivalent to solving a
numerical problem that computes a relevant permutation of \(Y\). To
solve it, we select the suitable permutations by minimizing a cost
function \(c(X,Y)\), as follows \begin{equation} \label{DE}
  \overline{\sigma} = \arg \inf_{\sigma \in \Sigma} c(X,Y\circ \sigma),
\end{equation} where \(\Sigma\) denotes the set of all possible
permutations.

We distinguish between two possibilities, as follows.

\begin{itemize}
\item
Assume first \(D_X=D_Y=D\).  In this case, the discrete problem is closely related to the Monge-Kantorovich problem, which has been extensively studied (see \cite{Brezis:2018} for the first result on the discrete transport problem). We employ non-Euclidean distance-based cost functions, typically defined with the discrepancy distance, say 
\begin{equation}
  c(X,Y) = \text{Tr}\big(M_k(X,Y)\big). 
\end{equation} 

\item  Assume next \(D_X \neq D_Y\). In this case, we propose to use a a cost functional based on the gradient, of the form
\begin{equation} 
\label{TC}
  c(X,Y) = \|\nabla Y_{k, \theta}(X) \|_2^2, \qquad \theta = K(X,X)^{-1} Y. 
\end{equation}
This approach is reminiscent of the \textit{travelling salesman problem}, and can be seen as a generalization of this classical problem.
\end{itemize}

\subsection{Multi-scale kernel
methods}\label{multi-scale-kernel-methods}

Let us consider the general extrapolation-interpolation formula
\eqref{FIT} and restrict attention to the extrapolation mode (so
\(\epsilon=0\)). Then the algorithmic complexity and memory cost of an
extrapolation with a kernel method consists of \(\mathcal{O}(N_X^2)\)
operations in order to store the Gram matrix, and \(\mathcal{O}(N_X^3)\)
operations in order to invert it.

These numbers of operations are unreasonable for large datasets, and we
thus propose a \textit{divide and conquer} strategy, which allows us to
handle large scale systems with kernel methods. Our approach consists in
interpreting \(N_Y\) in \eqref{FIT} as a number of computing units, each
one being defined by a \textit{centroid} denoted by \(y^n\), with
\(Y=(y^1,\cdots,y^{N_Y})\). This approach involves first the computation
of \(N_Y\) clusters that partition the distribution \(X\). By
considering hard partitioning (in which, by definition, each data point
is assigned to exactly one cluster), we denote by
\(l : \mathbb{R}^D \mapsto [1,\ldots,N_Y]\) the allocation function that
attaches a given point to a cluster.

The general construction is as follows. Consider the given cluster
centers \(Y\), as well as an allocation function \(l(\cdot)\). Choose
first a kernel at the coarser level \(k_0(\cdot,\cdot)\) and implement
\eqref{FIT}. This step generates the error term \[
  \epsilon(\cdot) = f(\cdot)-f_{k_0,\theta}(\cdot).
\] We then define the following extrapolation step, paving the way to an
easy implementation of a parallelized kernel algorithm based on \(N_Y\)
computation units: \begin{equation}\label{MS}
  f_{k_0,\cdot,k_{N},\theta}(\cdot) = f_{k_0,\theta}(\cdot) + \epsilon_{k_{l(\cdot)},\theta}(\cdot). 
\end{equation}

Here, each kernel \(k_n\) (for \(n=1,\cdot,N_Y\)) is used in the formula
\eqref{FIT} in extrapolation mode on each cluster, defined as a data
partition with an average size of \(\sim \frac{N_X}{N_Y}\) points.
Although the proposed method introduces an extra-parameter \(N_Y\), it
is quite appealing as it provides us with a rather natural way to design
a scalable, implementation-friendly, kernel method satisfying both the
reproductibility property and useful error estimates.

Observe that partitioning the dataset with \(N_Y\) clusters collecting
each exactly \(\frac{N_X}{N_Y}\) dataset points requires the use of
\textit{balanced clustering} methods. We are going to describe this
technique in Section \ref{balanced-clustering}, where we also specialize
these tools tothe application of multiscale methods to optimal
transportation.

In term of algorithmic complexity, the formula \eqref{FIT} requires
\(\mathcal{O}(N_X^3)\) operations to run in the extrapolation mode. The
constant depends on the time needed to perform each kernel evaluation,
which usually scales linearly with the dimension \(D\); therefore, we
obtain \(D\ln(D)\) for convolution type kernels. This approach requires
\(\mathcal{O} (\frac{N_X^3}{N_Y^3})\) operations for each of the \(N_Y\)
computational units, plus \(\mathcal{O} (N_X^2 N_Y + N_Y^3)\) for the
rough cluster. In particular, we have to take care of clustering
algorithms that run in a time less than the two first in order to avoid
to reduce the overall performance.

As data are separated, the computational units should ideally be run in
a parallel environment for \(N_Y\) threads, but they can also be run
sequentially on an architecture parallelizing each computational unit,
or combining both in order to easily adapt to a given computational
environment.

Of course, this approach has its own cost and somewhat reduces the
accuracy: the main contribution to the error for a general kernel
extrapolation method is given by the distance \(d_k(\cdot,X)\). The
overall accuracy of this multi-scale approach is thus
\(d_{k_0}(\cdot,Y)d_{k_{l(\cdot)}}(\cdot,X_{l(.)})\), where
\(X_i \subset X\) denote the data points in the \(i\)-th cluster. These
terms are computationally tractable, however we are not aware of
theoretical estimates for most of the available kernels. Our
experimental observations show a nonlinear increase in the dimension
\(D\), but a decrease at rates depending on the kernel \(k\), and the
inverse of the size of the dataset, that is \(N_X\), in accordance with
theoretical convergence rates available for simple situations. In other
words, with the proposed multi-scale approach, We trade accuracy for
computational time, the loss of accuracy depending on the geometry of
the problem at hand, the choice of kernels, and the parameter \(N_Y\).

From the engineering point of view, other quite interesting approaches
can be designed from this first construction: the first kernel can be
chosen as simple regressors, for instance a polynomial basis in order to
accurately match the moments of the distributions under consideration.
Other layers in a tree-like model can be designed to adapt to very
large-scale datasets. Using overlapping clusters between layers is also
an attractive strategy: one then seeks to use a first specialized kernel
that model some feature or interaction, and the error is passed to
another kernel that is better adapted to capture other features, acting
as filters. An example of such a construction for image recognition is
the use of a first convolution kernel, to match a pattern, the error
begin next treated by a rotation invariant kernel (for instance).

More generally, the strategy proposed here suggests that we could design
kernel methods by using oriented graphs

\section{General purpose algorithms}\label{general-purpose-algorithms}

\subsection{Greedy search algorithms}\label{greedy-search-algorithms}

We describe first several general purpose algorithms, on which we will
rely in order to implement the strategy outlined so far The following
greedy search algorithms provide us with a family of efficient and
versatile algorithms for the approximation of the following (quite
broad) class of problems: \begin{equation}
\label{GS}
\inf_{Y \subset X} D\big(Y,X\big), \quad D(Y,X) = \sum_n d(Y,x^n),
\end{equation} in which \(D(Y,X)\) denotes a (user defined) distance
measure between sets. We consider here a greedy search algorithm, that
is, \begin{equation}
\label{GS=2}
Y^{n+1} = Y^n \cup \arg \sup_{x \in X} d\big(Y^n,x\big). 
\end{equation} We will present some useful examples of this algorithm in
the context of kernel interpolation. For the literature on greedy
approximations, we refer the reader to
\cite{Temlyakov:2011,WaSaHa:2023,García:2024}.

This approach leads us to Algorithm \ref{cl0}, as stated, where the
parameter \(M\) (taken to be \(1\) by default) is introduced as a simple
optimization mechanism i order to trade accuracy versus time. Other
strategies can be implemented, more adapted to the problem at hand. Such
algorithms usually rely on optimization techniques within the main loop
for faster evaluation of \(d\big(Y^n, \cdot\big)\), using
pre-computations or \(d\big(Y^{n-1}, \cdot\big)\).

\begin{algorithm}
\caption{Greedy Search Algorithm} \label{cl0}
\LinesNumbered
\SetAlgoLined
\SetKwBlock{Begin}{}{end}
    \KwIn{a training set $X$, a distance measure $d(\cdot,\cdot)$, two integers: $1 \le N_Y \le N_X$, where $N_Y$ denotes the number of clusters, $M$ is an optional \textit{batch} number (taken to $1$ by default), and $Y^0 \subset X$ is  an optional set of initial points.}
    \KwOut{A set of indices $\sigma : [1...N_Y] \mapsto [1...N_X]$, defining $X\circ\sigma \subset X$ of size $N_Y$ as approximate clusters.}
            for $n=0,\ldots,N_Y/M$ 
            \Begin{
             find a new numbering, say $X^n$, according to the decreasing order of $d\big(Y^n, x^p\big)$, $p=1,\ldots,N_X$ \\
              $Y^{n+1}=Y^{n}\cup X[1,\ldots,M]$.
    }            
\end{algorithm}

\subsection{Permutation algorithms}\label{permutation-algorithms}

\subsubsection{Linear sum assignment
problems}\label{linear-sum-assignment-problems}

The so-called LSAP problem (Linear Sum Assignment Problem) is a
cornerstone of combinatorial optimization, with wide-ranging
applications across academia and industry. The problem has been
extensively studied and is well
documented\footnote{ \url{https://en.wikipedia.org/wiki/Assignment_problem}}.
It can be described as follows: \begin{equation}
\label{LSAP}
  \overline{\sigma} = \arg \inf_{\sigma \in \Sigma} \sum_{i=1}^M c(i,\sigma^i), 
\end{equation} where \(\Sigma\) is the set of all injective renumberings
\(\sigma:[1,\ldots,M] \mapsto [1,\ldots, N]\) with \(M \le N\). Of
course, when \(M=N\), this is nothing but the set of permutations. These
algorithms \textit{put together} two distributions according to a
similarity criterion given by a general, rectangular \textit{cost} or
\textit{affinity} matrix \(C \in \mathbb{R}^{M,N}\).

\subsubsection{Discrete descent
algorithms}\label{discrete-descent-algorithms}

Let us explore a generalized version of this linear sum assignment
problem, specifically the following optimization problem:
\begin{equation}\label{GLSAP}
\bar{\sigma} = \arg \inf_{\sigma \in \Sigma} C(\sigma), 
\end{equation} where \(\Sigma\) denotes the set of all possible
permutations, and \(C(\sigma)\) is a general cost functional. A special
case of this problem is indeed the LSAP above, which corresponds to a
linear cost function \(C(\sigma) = \sum_i c(i, \sigma^i)\). However, in
the present paper, we do consider other forms of cost functionals which
may be better tailored for our applications.

For problems where the permutation gain, defined as
\(s(i,j,\sigma) = C(\sigma)-C(\sigma_{i,j})\) can be efficiently
computed, we use a simple descent algorithm. Here \(\sigma_{i,j}\)
represents the permutation obtained by swapping the indices \(\sigma^i\)
and \(\sigma^j\). This approach leverages the fact that any permutation
\(\sigma\) can be decomposed into a sequence of elementary two-element
swaps. For instance, for LSAP problems\eqref{LSAP}, the permutation gain
is given by
\(s(i,j,\sigma) = c(i,\sigma^i)+c(j,\sigma^j)-c(i,\sigma^j)- c(j,\sigma^i)\).
In its simplest form, this algorithm is summarized in Algorithm
\ref{PBA}. Observe that there exist more performing algorithms than this
discrete descent approach, but adapted to specific situations; see
\cite{TH:2019}. Adapting them to our purpose is a research in progress.

\begin{algorithm}
\caption{Permutation-Based Descent Algorithm} 
\label{PBA}
\begin{algorithmic}[1]
\REQUIRE A permutation-gain function $s(i,j,\sigma)$, where $\sigma:[1,\ldots,N_J] \mapsto [1,\ldots,N_I]$ is any injective mapping (that is, a permutation if $N_I=N_J$).
\vskip.15cm
\ENSURE An injective mapping $\sigma:[1,\ldots,N_J] \mapsto [1,\ldots,N_I]$ achieving a local minima to \eqref{GLSAP}.
\WHILE{FLAG = True} 
\STATE FLAG $\leftarrow$ False
\FOR{$i=1 \ldots, N_I$,$j=1 \ldots, N_J$} 
\IF{$s(i,j,\sigma) < 0$}
\STATE $\sigma \leftarrow \sigma^{i,j}$
\STATE FLAG $\leftarrow$ True
\ENDIF
\ENDFOR
\ENDWHILE
\end{algorithmic}
\end{algorithm}

The \texttt{For} loop in this algorithm can be further adapted to
specific scenarios or strategies, for instance a simple adaptation to
symmetric permutation gain function \(s(i,j,\sigma)\) (as met with the
LSAP problem) is given as follows:
\texttt{for\ i\ in\ {[}1,\ N{]},\ for\ j\ \textgreater{}\ i} in order to
minimize unnecessary computation.

While these algorithms generally produce sub-optimal solutions for
non-convex problems, they are robust and tend to converge in finite
time, i.e.~within a few iterations of the main loop, and a careful
choice of the initial permutation \(\sigma\) can escape sub-optimality.
They are especially useful as auxiliary methods in place of more complex
global optimization techniques or for providing an initial solution to a
problem. Additionally, they are advantageous in finding a local minimum
that remains \emph{close} to the original ordering, preserving the
inherent structure of the input data.

However, these algorithms have some limitations. Depending on the
formulation of the permutation gain function \(s(i,j,\sigma)\), it is
possible to parallelize the \texttt{For} loop. Nonetheless,
parallelization often requires careful management of concurrent
read/write access to memory, which can complicate implementation and
alter performances. Furthermore, theoretical bounds on the algorithm
complexity are typically not available, making performances prediction
difficult in all cases.

\subsection{Explicit descent
algorithms}\label{explicit-descent-algorithms}

We now present a generalized gradient-based algorithm, specifically
designed for the minimization of functionals of the form
\(\inf_X J(X)\), where the gradient \(\nabla J(X)\) is locally convex
and is explicitly known. In this scenario, we can apply the simplest
form of a gradient descent scheme, often referred to as an Euler-type
method. In its continuous form, it is written as
\(\frac{d}{dt} X(t) = -\nabla J(X(t))\), \(X(0)=X^0\) and, in its
numerical formn, it takes the form \begin{equation}
\label{DA}
  X^{n+1} = X^{n} - \lambda^n\nabla J(X^{n}).
\end{equation} The term \(\lambda^n\) is known as the
\textit{learning rate}. In this situation, as the gradient of the
functional is explicitly given, we can compute sharp bounds over
\(\lambda^n\), allowing to apply root-finding methods, such as the Brent
algorithm\footnote{For instance \href{https://en.wikipedia.org/wiki/Brent's_method}{this wikipedia's page}},
and efficiently locate the minimum while avoiding \textsl{instability}
issues often met with Euler schemes.

\begin{algorithm}
\caption{Descent Algorithm With Explicit Gradient} 
\label{DA=algo}
\begin{algorithmic}[1]
\REQUIRE A function $J(\cdot)$, its gradient $\nabla J(\cdot)$, a first iterate $X^0$, tolerance $\epsilon >0$ or number of maximum iterations $N$.
\vskip.15cm
\ENSURE A solution $X$ achieving a local minimum of $J(X)$.

\WHILE{$\|\nabla J(X^n)\| > \epsilon$ or $n < N$} 
\STATE compute $\lambda^{n+1} = \arg \inf_{\lambda} \|\nabla J(X^\lambda)\|$, where $X^{\lambda}:= X^{n} - \lambda \nabla J(X^{n})$ with a root-finding algorithm.
\STATE $ X^{n+1} = X^{n} - \lambda^{n+1}\nabla J(X^{n})$
\ENDWHILE
\end{algorithmic}
\end{algorithm}

\section{Clustering algorithms}\label{clustering-algorithms}

\subsection{Proposed strategy}\label{proposed-strategy}

Clustering is regarded as a key exploratory data analysis task, which
aims at collecting objects into clusters on the basis of a similarity
criterion. In the context of kernel methods, clustering serves an
additional crucial purpose: reducing the computational burden of kernel
operations. Indeed, having a look to the error formula \eqref{dk}, the
best choice to determine the centers \(Y\) at a coarse level in order to
approximate a given kernel-induced functional space
\(\mathcal{H}_k^{X}\) is given by the following minimization problem:
\begin{equation} \label{SDS}
    \arg \inf_{Y \subset \RR^{N_Y,D}} d_k\big(Y,X\big)^2.
\end{equation} We call such sequences \textit{sharp discrepancy}
sequences, since they define the best possible meshes for kernel
methods, namely a discrete or continuous distribution \(X\). They can be
explicitly computed in sufficiently simple cases, with Fourier-based
kernels; see the discussion in \cite{PLF-JMM-estimate}; however we rely
on a numerical approximation of \eqref{SDS} for most of our kernels. We
can easily design descent-type methods in order to solve \eqref{SDS}.
Unfortunately, such a direct approach can be troublesome: descent
algorithms are difficult to scale up to large datasets and dimensions,
if one chooses to avoid stochastic descent approaches, and one can be
trapped at a local minimum since we are dealing with a strongly non
convex problem. Note that the minimization problem \eqref{SDS} considers
the full discrepancy functional \(d_k\big(Y,X\big)^2\), so we are
considering a different clustering method that what is referred in the
literature as kernel k-means, a natural extension of k-means algorithms
to \eqref{SDS}.

Hence, we prefer to rely upon faster, but less accurate, algorithms and
fasten the computation of such sequences. We mainly use three
algorithms, presented in the order of their numerical burden, fastest to
heaviest, the last one being a descent algorithm to solve \eqref{SDS}.
We generally also combine them in this order: the output of an algorithm
becomes the input of the following one.

The philosophy behind our strategy is to seek first approximate clusters
in the form of subset \(Y \subset X\) of the training set, where we can
design fast combinatorial algorithms, which do benefit from a
preliminary computation. Then, eventually, we seek for clusters outside
the training set, solving \eqref{SDS} with descent algorithms. In the
following discussion, \(N_Y\) refers to the number of clusters, and
\(N_X\) to the size of the training set.

\subsection{Greedy clustering methods}\label{greedy-clustering-methods}

Our first algorithm is based on considering, instead of \eqref{SDS}, the
following problem: \begin{equation} 
\label{SDSD}
    \inf_{Y \subset X} d_k\big(Y,X\big),
\end{equation} tackling it with the Greedy-Search Algorithm \ref{cl0}.
To this aim, elementary computations leads to define the following
functional to use with Algorithm \ref{cl0}

\[
    d(Y,\cdot) = \frac{1}{(N_Y+1)^2}\Big(2 \sum_{m=1}^{N_Y} k(y^m,\cdot)+k(\cdot,\cdot)\Big) - \frac{2}{(N_Y+1)N_X}\sum_{n=1} ^{N_Y} k(x^n,\cdot).
\] The terms \(\sum_{n=1}^{N_x}  k(x^n,x^m)\) can be precomputed,
requiring \(\mathcal{O}(N_X^2)\) calls to the kernel, with a memory
footprint of \(\mathcal{O}(N_X)\). At each step \(n\) of Algorithm
\ref{cl0}, the functional \(d(Y^n,\cdot)\) requires
\((n \times N_X -n)\) calls to the kernel function, which can be run in
parallel, or a look-up to a pre=computed matrix of size \(N_X^2\) if the
problem is sufficiently small in order to fit into the available memory.
The leading complexity terms are \(\mathcal{O}(N_Y^2 N_X + N_X^2)\).
This algorithm can handle medium-to-large datasets, depending on
settings. In this paper, we consider a version that necessitate to
precompute the full Gram matrix, requiring \(\mathcal{O}(N_X^2)\) calls
to the kernel, hence is adapted to medium datasets. However, kernel
k-means algorithms, which are version of k-means algorithms working with
the discrepancy distance (instead of the Euclidean one), should be able
to handle large scale datasets, but compute clusters with lower quality.

Solving the \textit{sharp discrepancy} problem \eqref{SDS} is useful for
mesh generation techniques, when one is more interested in a relevant
approximation of the whole space \(\mathcal{H}_k^{X}\). It relates to
unsupervised learning. However, if one considers a single function \(f\)
instead of the whole functional (as in the interpolation formula
\eqref{FIT}), the context is supervised learning, and can also be
tackled with a greedy search algorithm. To that aim, we define the
following distance for any \(1 \le p\le \infty\): \begin{equation} 
\label{DF}
 D(Y,X) = \sum_{n=1}^{N_X} \big|f(x^n) - f_{k, \theta_Y}(x^n)\big|^p, \quad \theta_Y:=K(Y,Y)^{-1}f(Y),
\end{equation} where \(f\) is a given, vector-valued function. Solving
\(\inf_{Y \subset X} D\big(Y,X\big)\) accounts to selecting from the
distribution \(X\) those points that best represent \(f\). This relates
to kernel \textit{adaptive mesh}, or kernel \textit{control variate}
techniques, which are useful for numerical simulations as well as
statistics.

Considering \eqref{DF}, we use the greedy algorithm in combination with
block matrix inversion techniques to fasten the computation of
\(K(Y^n, Y^n)^{-1}\) from \(K(Y^{n-1}, Y^{n-1})^{-1}\). This small
optimization maintains the overall algorithmic complexity at a rate
\(\mathcal{O}(N_Y^3 + N_X N_Y^2)\) with minimal memory requirements,
making it suitable to approximate a large number of clusters on large
datasets.

\subsection{Subset clustering methods}\label{subset-clustering-methods}

Our next clustering algorithm allows us to refine the previous
approximation of \eqref{SDSD}, using the Discrete Permutation Algorithm
\ref{GLSAP}. Let us define a permutation
\(\sigma:[0\ldots N_X] \mapsto [0\ldots N_X]\), and, according to
\eqref{dkerr}, the permutation following gain function \[
  s(i,j,\sigma) := K_{XY}^{\sigma}-K_{XY}^{\sigma_{i,j}}+K_{YY}^{\sigma}-K_{YY}^{\sigma_{i,j}},
\] where \(\sigma_{i,j}\) is the permutation obtained from swapping
\(i\) and \(j\), the vectors \(K_{XY}^{\sigma},K_{YY}^{\sigma}\) are
pre-computed from contracted matrices \[
    K_{YY}^{\sigma} := \frac{1}{N_Y^2} \sum_{n,m}^{N_Y,N_Y} k(x^{\sigma^n},x^{\sigma^m}), \quad K_{XY}^{\sigma}:= - \frac{2}{N_X N_Y} \sum_{n,m}^{N_X,N_Y} k(x^{n},x^{\sigma^m}).
\] This defines clusters as \(Y = X\circ \sigma(1 \ldots N_Y)\), and we
can set \(s(i,j,\sigma)=0\), if \(i\ge N_Y\) and \(j\ge N_Y\), or
\(i\le N_Y\) and \(j\le N_Y\). We thus can use the permutation algorithm
\eqref{PBA}, adapting the for loop as follows
\texttt{for\ i\ in\ {[}1,\ N\_Y{[},\ for\ j\ in\ {[}N\_Y,\ N\_X{[}}.
This allows to quickly find a sub-optimal solution to \eqref{SDSD},
defining an approximation of \(Y\). To avoid sub-optimal solutions in
this descent algorithm as much as possible, it is beneficial to start
with a pertinent choice of the input set \(Y\). A suitable first choice
can be obtained using the greedy algorithm described in the previous
section.

While we do not provide an explicit estimate of the algorithm
complexity, empirical observations suggest that it provides us with
polynomial complexity, approximately \(\mathcal{O}(N_X^2 N_Y)\), and can
be run parallel. However this algoritmic complexity is an empeachment to
large dataset, where kernel k-means algorithms should be used instead.

\subsection{Sharp discrepancy
sequences}\label{sharp-discrepancy-sequences}

Finally, we discuss clustering sequences that aim to approximate
\textit{sharp discrepancy} sequences, which are solutions to
\eqref{SDS}. For this problem, one can use the general descent Algorithm
\ref{DA=algo}, with \(J(Y)=d_k(X,Y)\), for which the expression
\eqref{nabla}, applied to \(d_k\big(Y,X\big)^2\), provides an exact
formula of the gradient of this functional's \(\nabla J(Y)\), to
approximate a local minimum solution. To avoid local minima, it is
advisable to consider Algorithm \ref{DA=algo} with first guess \(X^0\)
given by the subset algorithms described previously.

Two remarks are in order:

\begin{itemize}
\item
  Such an algorithm computes very good clustering sequences, but is
  computationally heavy and might be difficult to apply to large scale
  dataset and dimensions. The previous combinatorial algorithms provide
  faster alternatives, at the cost of less qualitative clusters
  retrieval.
\item
  For a number of non-smooth kernels \(k\), the discrepancy functional
  is \textbf{concave} almost everywhere. Continuous descent algorithms
  as Algorithm \ref{DA=algo} are not adapted to this situation, and
  combinatorial algorithms should be preferred while looking for a
  global minimum in this context.
\end{itemize}

\subsection{Balanced clustering}\label{balanced-clustering}

\subsubsection{Balanced clustering for supervised
learning}\label{balanced-clustering-for-supervised-learning}

For our multiscale kernel method, it is paramount to define clusters of
comparable size. There exist several algorithms capable to handle
balanced clusters, see for instance \cite{MaFr:2014}. In our context, we
work with cluster algorithms based upon an induced distance
\(d(\cdot,\cdot)\), that is the euclidean distance for k-means and the
kernel discrepancy for kernel-based clustering algorithms. Hence we
designed a general approach to balanced clustering, that takes as input
a distance matrix
\(D \in \mathbb{R}^{N_Y,N_X}=\Big(d(y^i,x^j)\Big)_{i,j}\), where the
centroids \(y^i\) are given by another algorithm. Our balanced cluster
method amounts to solve the following discrete optimal transport
problem, where \(\%\) holds for modulo \begin{equation}\label{balanced}
  \inf_{\sigma \in \Sigma} \sum_{n=1}^{N_X} d(y^{(\sigma^n\ \% \ N_Y)},x^n).
\end{equation} Such a discrete minimization problem can be solved with
the Discrete Descent Algorithm \ref{PBA}, together with the permutation
gain function \[
  \sigma(i,j,\sigma)=D(i,\sigma^{i \% N_Y})-D(j,\sigma^{j \% N_Y}) - D(i,\sigma^{j \% N_Y}) +D(j,\sigma^{i \% N_Y}).
\] This method is balanced, by construction, as it attaches each point
\(x^n\) to the cluster \(y^{\sigma^n \ \% \ N_Y}\), and it also can
provides us with an allocation function, defined by
\(l(\cdot) = l^{\arg \inf_n d(x^n,\cdot)) \ \% \ N_Y}\). Observe that
this algorithm can also treats \(N_Y\) random points among \(X\),
together with a user-defined distance, such as the Euclidean or
discrepancy distances. With such a naive setting, this method computes
balanced clusters and, moreover, as highlighted in our numerical
experiments, is numerically efficient and provide us with clusters
having good quality.

\subsubsection{Balanced clustering for optimal
transport}\label{balanced-clustering-for-optimal-transport}

We precise in this paragraph the clustering method used to approximate
the optimal transportation \eqref{ED} with \(X,Y \in \mathbb{R}^{N,D}\)
for large dataset size \(N\). Here, a regressor \(Y_{k,\theta}(\cdot)\)
computed by the multiscale formula \eqref{MS} defines an exact transport
of \(X\) into \(Y\), that is satisfying \(Y_{k,\theta, \#}X=Y\), as
would do \(Y\circ \sigma_Y\), where
\(\sigma_Y : [0,\ldots,N]\mapsto [0,\ldots,N]\) is any permutation.

Multiscale method for optimal transport includes a hidden difficulty :
let \(M\) the clusters number, and denote two centroids set
\(C_X,C_Y \in \mathbb{R}^{M,D}\), one for each distribution, computed by
an external clustering method. We not only need to find an association
\(\sigma_C : [0,\ldots,M] \mapsto [0,\ldots,M]\) in order to map each
cluster into another, but we also need clusters of equal sizes, as
required by the extrapolation formula \eqref{FIT}. Thus it is mandatory
to deal with an adapted balanced clustering method.

To that aim, let us consider the following problem, where
\(d(\cdot,\cdot)\) is a general distance function, for instance the
Euclidean one or the discrepancy one \(d_{k_0}(\cdot,\cdot)\), where
\(k_0\) is the coarser kernel used in the multiscale formula \eqref{MS},
\(\%\) holds for modulo, and \(\Sigma\) is the set of all permutations
of \([0,\ldots,N]\):

\begin{equation} \label{KM}
  \overline{\sigma_X},\overline{\sigma_Y} = \arg \inf_{\sigma_X,\sigma_Y \in \Sigma} \sum_{n} D_X(n,\sigma_X^{n \% M}) + \sum_{n} D_Y(n,\sigma_Y^{n \% M})+\sum_n C(\sigma_X^{n \% M},\sigma_Y^{n \% M} )
\end{equation} where we denoted the user defined gram matrix
\(D_X = d(X,C_X)\), \(D_Y = d(Y,C_Y)\), \(C = D(C_X,C_Y)\) having size
resp. \((N,M)\),\((N,M)\),\((M,M)\).

Once computed or approximated, this problem defines

\begin{itemize}
\item Two affectation functions $l_X(\cdot) = \overline{\sigma_X}^{\arg \inf_n d(X^n,\cdot)\% M}$ (resp. $l_Y(\cdot) = \overline{\sigma_Y}^{\arg \inf_n d_{k_0}(Y^n,\cdot)\% M}$), attaching any point to the centroid $C_X^{l_X(\cdot)}$ (resp. $C_Y^{l_Y(\cdot)}$).
\item The association $C_X^{\overline{\sigma_X}^n \ \% M} \mapsto C_X^{\overline{\sigma_Y}^n \ \% M}$
\item The clusters are balanced by construction.
\end{itemize}

Once this problem solved, we define the approximation of the optimal
transport map as for the supervised case \eqref{MS}:

\begin{equation}\label{MSOT}
  Y_{k_0,\cdot,k_{N},\theta}(\cdot) = Y_{k_0,\theta}(\cdot) + \epsilon_{Y_{l(\cdot)},\theta}(\cdot), \quad \epsilon(\cdot) = Y - Y_{k_0,\theta}(\cdot).
\end{equation}

As we split the dataset, this method consists in an approximation of the
optimal transport map, but defining an exact transportation map, that is
invertible if restricted to each clusters.

\section{A framework for dealing with two
distributions}\label{a-framework-for-dealing-with-two-distributions}

\label{A-framework-for-dealing-with-two-distributions}

\subsection{Encoder-decoders and sampling
algorithms}\label{encoder-decoders-and-sampling-algorithms}

We now turn to algorithms related to the distributional framework
described in Section \ref{distrib}, which have proven useful for our
applications. The first algorithm we present implements the equation
\eqref{DE}, where, given matrices \(X \in \mathbb{R}^{N, D_X}\) and
\(Y \in \mathbb{R}^{N, D_Y}\), we seek to define continuous, invertible
maps that push \(X\) into \(Y\). As mentioned in \eqref{ED}, we consider
two cases:

\begin{itemize}
\tightlist
\item
  If \(D_X = D_Y\), the problem is a Monge-Kantorovich transport problem
  and can be solved using the LSAP algorithm.
\item
  If \(D_X \neq D_Y\), we solve the following optimization problem:
  \begin{equation}  
  \label{perm2}
  \sigma_* =   \arg \inf_{\sigma \in \Sigma } \| \nabla (Y\circ \sigma)_{k, \theta}(X) \|^2_{\ell^2} = \arg \inf_{\sigma \in \Sigma } \langle -\Delta_k , (Y\circ \sigma) (Y\circ \sigma)^T \big)\rangle,
  \end{equation} where \(\Delta_k\) is the discrete Laplacian operator
  defined in \eqref{delta} and \(\langle \cdot, \cdot \rangle\) is the
  Frobenius inner product. We use this observation to design the Kernel
  Sampling Algorithm \ref{alg1}, as follows.
\end{itemize}

\begin{algorithm}
\caption{Kernel Sampling Algorithm} \label{alg1}
\begin{algorithmic}[1]
\item[\textbf{Settings:}] A kernel $k$
\REQUIRE data $X \in \mathbb{R}^{N,D_X}$, $Y \in \mathbb{R}^{N,D_Y}$, with $x^i\neq x^j$, $y^i\neq y^j$,  $\forall \ i \neq j$
\ENSURE A regressor $Y_{k, \theta}(\cdot)$ satisfying $Y_{k, \theta}(X) = Y$, modeling an invertible push-forward map.
\IF{$D_X=D_Y$}
  \STATE Compute a permutation $\sigma$ using LSAP \eqref{LSAP} with cost function $c(i,j)=d_k(x^i,y^j)$
\ELSE
  \STATE Perform precomputations in order to provide a fast evaluation of the permutation gain function
  $$s(i,j,\sigma)=<\Delta_k,(Y\circ\sigma) (Y\circ\sigma)^T> - <\Delta_k,(Y\circ\sigma_{i,j}) (Y\circ\sigma_{i,j})^T>.$$
  \STATE Compute a permutation $\sigma$ using Algorithm \ref{PBA}.
\ENDIF
\STATE Return the regressor $(Y\circ \sigma)_{k, \theta}(\cdot)$, where $\theta$ is computed with \eqref{ED}.
\end{algorithmic}
\end{algorithm}

This algorithm is primarily used to generate samples from a distribution
\(\mathbb{Y}\), known only through a sample
\(Y \in \mathbb{R}^{N, D_Y}\). One can use a known distribution
\(\mathbb{X}\) (e.g., a normal distribution), draw a sample
\(X \in \mathbb{R}^{N, D_X}\), and compute the regressor
\((Y \circ \sigma)_{k, \theta}(\cdot)\). The map
\(z \mapsto (Y \circ \sigma)_k(z, \theta)\), where \(z\) is drawn from
\(\mathbb{X}\), provides a generator for producing new samples of
\(\mathbb{Y}\). This process resembles generative methods, such as those
used in GAN architectures. In this context, the Artificial Intelligence
community uses specific terminology:

\begin{itemize}
    \item \textbf{Latent space:} The space $\mathbb{R}^{D_X}$ is referred to as the latent space. 
    \item \textbf{Decoding:} The map $x \mapsto (Y \circ \sigma)_{k, \theta}(x)$ is called a decoder.
    \item \textbf{Encoding:} The map $y \mapsto (X \circ \sigma)_{k, \theta}(y)$, with $\theta = K(Y \circ \sigma, Y \circ \sigma)^{-1} X$, is called an encoder.
\end{itemize}

The latent space is typically thought of as a lower-dimensional space
that captures the essential features of the data. This provides a
compact, informative representation of the original data \(Y\), allowing
for efficient encoding and decoding of information while leveraging the
structural properties of the data.

The choice of the latent space \(\mathbb{R}^{D_X}\) is important for
applications, though it is not well documented. When using a standard
normal distribution and a small latent dimension (e.g., \(D_X = 1\)),
the resulting generator \(x \mapsto (Y \circ \sigma)_{k, \theta}(x)\)
produces samples that closely resemble the original sample \(Y\), which
can pass standard statistical tests such as the Kolmogorov-Smirnov test.
Using larger dimensions (e.g., \(D_X \gg 1\)) results in greater
variability in the samples, but may cause the model to fail statistical
tests as the dimension increases. A possible explanation for this
phenomenon is that the average kernel distance \(d_k(x^i, x^j)\) between
two draws from standard distributions, such as the normal distribution,
increases with the latent space dimension \(D_X\).

One of the major advantage of this kernel generative methods is to
produce continuous generators, that can exactly reconstruct the original
variate \(Y\) if needed, as the reproductibility property \eqref{REPRO}
ensures \((Y\circ \sigma)_{k, \theta}(X) = Y\).

\subsection{Conditional distribution sampling
model}\label{conditional-distribution-sampling-model}

The generation of conditioned random variables is an essential part of
the encoder-decoder framework. Given two variates
\(X \in \mathbb{R}^{N, D_X}\) and \(Y \in \mathbb{R}^{N, D_Y}\) from two
distributions, our goal is to provide a generator modeling the
conditioned distribution \(Y|X=x\).

Considering a latent variable having form
\(\eta:=(\eta_x,\eta_y) \in \mathbb{R}^{N,D_{\eta}}\), with
\(D_{\eta} = D_{\eta_x}+D_{\eta_y}\), the algorithm \eqref{alg1} allows
to determine this two mappings:

\begin{itemize}
    \item An encoder targeting the latent distribution $\eta_x$, $x \mapsto (\eta_x\circ \sigma_X)_{k,\theta}(x)$.
    \item A decoder targeting the joint distribution $(X,Y)$, $\eta \mapsto ([X,Y] \circ \sigma_{XY})_{k,\theta}(\eta)$.
\end{itemize}

Our conditioned model consists to consider the following generator as a
model for the distribution \(Y|X=x\): \begin{equation} 
\label{condalg}
  \eta_y \mapsto (Y\circ \sigma)_{k,\theta}([\eta_x,\eta_y]), \quad \eta_x = (\eta_x\circ \sigma)_{k,\theta}(x),
\end{equation} where the notation
\(\eta \mapsto (Y \circ \sigma)_{k,\theta}(\eta)\) refers to the
\(Y\)-component of the joint distribution. Some remarks concerning this
construction :

\begin{itemize}
    \item In some situations, there is no need to define the encoder $\eta_x\circ \sigma_X$, and $x$ is then considered as a latent variable.
    \item We can extend this approach to model more elaborated conditioning, as multiple conditioning $Y|X \in \{x^1,\ldots,x^n\}$, or more generally $Y|\mathbb{Z}$, $\mathbb{Z}$ being a distribution supported in $\mathbb{R}^{D_X}$.
\end{itemize}

\subsection{Transition probability
algorithms}\label{transition-probability-algorithms}

\paragraph*{Motivation}

The goals of these algorithms are to compute conditional expectations in
arbitrary dimensions. This is a fundamental problem having many
applications in stochastic analysis and ranking, where we often need to
estimate conditional expectation of a function \(g(Y)\), knowing that a
process \(t\mapsto X_t\) was in in the states \(X = X_{t_1}\),
\(Y=X_{t_2}\) at times \(t_1<t_2\), where \(g(\cdot)\) is an arbitrary,
vector-valued function: \begin{equation}
g(Y \mid X) := \mathbb{E}[g(Y) \mid X], \quad X,Y \in \mathbb{R}^{N , D}.
\end{equation} The challenge is to compute these conditional
expectations efficiently, particularly in high-dimensional spaces. To
adress this problem, we first make a structural hypothesis, that is
\begin{equation}
  t \mapsto X_t \text{ is a martingale process, for } t_1 \le t \le t_2
\end{equation} Under this hypothesis, we can write w.l.o.g
\(\mathbb{E}[X]=\mathbb{E}[Y]=0\) and approximate \begin{equation}
g(Y \mid X) = \Pi(Y \mid X)g(Y),
\end{equation} where \(\Pi(Y \mid X)\) is a doubly stochastic matrix,
that is \(\Pi:= (\pi_{i,j})_{i,j=1}^N\), with
\(\sum_i \pi_{i,j} = \sum_j \pi_{i,j} = 1\), for any \(i,j\), that we
write shortly as \(\Pi \ 1 = 1\), \(1 \ \Pi  = 1\) in the sequel. In the
context of martingale processes, these matrix are positive, that is
\(\pi_{i,j} \in [0,1]\), and are interpreted as the probability for the
process \(Y_{t}\) to \textit{jump} from the state \(X^i\) to the state
\(Y^j\).

We present here two algorithms, which should be considered as an
alternative to existing algorithms for the computation of transition
probability matrices. To our knowledge, there exist only two alternative
techniques. On the one hand, relying on an entropy-based regularization,
the Sinkhorn method \cite{Cuturi:2013} is a scalable and parallelizable
algorithm for solving optimal transport problems. On the other hand,
Neural networks also provide deep learning techniques; see
\cite{HG:2020} and the references cited therein.

Our first algorithm is the following version of the celebrated Nadaraya
Watson density estimator : let \(k\) a kernel, then this method
approximates the transition probability matrix as

\begin{equation} \label{NW}
  \Pi(Y \mid X) \sim \text{IPF}\Big(K(X,X) \cdot K(Y,Y)\Big)
\end{equation} where IPF stands for the iterative proportional fitting
algorithm\footnote{see, for instance, the wikipedia page https://en.wikipedia.org/wiki/Iterative\_proportional\_fitting and the references therein.}
Here, \(\cdot\) holds for the inner product of the two Gram matrices.

The second is a version of an algorithm introduced in
\cite{LM-CRAS}-\cite{PLF-JMM:2023}. This algorithm considers the
following minimization problem: \begin{equation} \label{MinPi}
  \inf_{\Pi \in \mathbb{R}^{N,N}} J(\Pi), \quad J(\Pi) = \| X - \Pi Y \|_{\ell^2}^2,\quad \text{constr. }  \Pi = \text{bi-stochastic}
\end{equation} Considering a continuous path \(t\mapsto \Pi_t\),
starting with \(\Pi_0 = I_d^N\), this descent algorithm produces the
following continuous-in-time scheme
\(\frac{d}{dt} \Pi = (X - \Pi Y)Y^T\), that we can integrate at first
order as \(\Pi_t = I_d^N + t(X - Y)Y^T\). Note that \(\Pi_{t}\) is
bi-stochastic, as \(\Pi_t \ 1 = 1 + t(X -  Y)Y^T \ 1 = 1\), since
\(\mathbb{E}[Y]=Y^T \ 1=0\). Since our target is to find
\(\bar{t} = \arg \inf_t J(\Pi_t)\), elementary computations provide the
minimum as the solution of the following iteration scheme
\begin{equation}\label{PI1}
\Pi_{\bar{t}}^{n+1} = \Pi_{\bar{t}}^{n} + \bar{t}^n(X -  Y^{n,T})Y^{n,T}, \quad Y^{n+1}=\Pi_{\bar{t}}^{n+1}Y, \quad t^n = \frac{<X-Y^n,(X -  Y^n)Y^{n,T}Y^n>}{|(X -  Y^n)Y^{n,T}Y^n|^2}
\end{equation} The motivation of this algorithm is to look to the
transition probability matrix as a perturbation of the identity one,
thus we first consider the reindexing procedure given by the LSAP method
\eqref{LSAP} in the corresponding Algorithm \ref{sumalg}. It provides a
fast and precise alternative to approximate a stochastic matrix,
although the output matrix might not be positive.

\begin{algorithm}
\caption{$\Pi$ algorithm}
\begin{algorithmic}[1] \label{sumalg}
\item[\textbf{Settings:}] Kernel $k$, $X,Y \in \mathbb{R}^{N,D}$ two sets of points of equal length, tolerance $\epsilon >0$ or maximum iteration number $M$.
\REQUIRE $X=X(t_1)$, $Y=X(t_2)$, where $t_1 \le t_2$ and $t \mapsto X(t)$ is a martingale process
\ENSURE $\Pi(X, Y)$ a bi-stochastic approximation of the transition probability matrix $\mathbb{P}(y^i | x^j)$

\STATE $Y^0 = Y\circ \sigma$, $\sigma$ being computed with the LSAP algorithm \eqref{alg1}.
\STATE $\Pi^0 = I_N$, the identity matrix.
\WHILE{ $\|Y^n-Y^{n-1}\| > \epsilon$ or $l < M$} 
\STATE compute $\Pi^{n+1},Y^{n+1}$ using \eqref{PI1}
\ENDWHILE
\STATE return $\Pi^{n+1}\circ \sigma^{-1}$. 
\end{algorithmic}
\end{algorithm}

\section{Numerical results}\label{numerical-results}

\subsection{Examples of clustering}\label{examples-of-clustering}

First, we illustrate some properties of our different clustering
methods, benchmarking them against Scikit implementation of the k-means
algorithm, which is a parallel, batch-based, highly optimized, and
scalable version of the standard k-means algorithm.''

We present our clustering methods for the standard Gaussian kernel
\(k(x,y) = \exp(|x-y|_2^2)\). Kernel methods typically require the use
of mapping functions to fit the data. For Gaussian kernels, we thus use
\(k\circ S\), where the map \(S\) is adapted to this kernel:
\(S=S_1\circ S_2\circ S_3\), where \(S_3(\cdot)\) maps the input data
\(X\) to the unit cube, \(S_2(\cdot)= \text{erf}^{-1}(\cdot)\)
(\(\text{erf}\) being the standard error function), and
\(S_1(\cdot) = \frac{\cdot}{\alpha}\) (a variance-type rescaling
\(\alpha=\frac{1}{N_XN_Y}\sum_{n,m} |x^n-y^m|_2^2\)).

\begin{figure}
\centering
\includegraphics{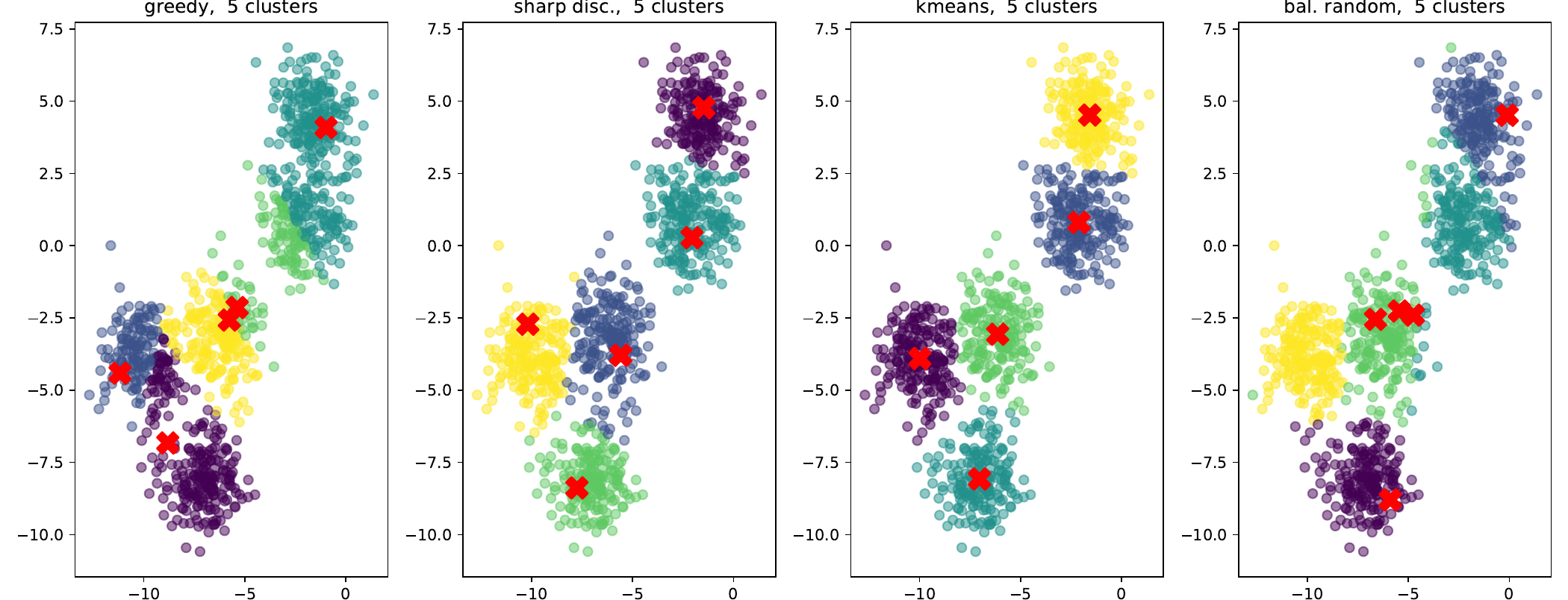}
\caption{Comparison of clustering methods}\label{fig:ccm}
\end{figure}

Our first illustration is qualitative, visually checking the quality of
the computed clusters \(Y\) for a given distribution \(X\) in two
dimensions. To this end, we generate a ``blob'', multi-modal
multivariate Gaussian with five modes. We then apply one of several
clustering methods and present the results in Figure \ref{fig:ccm}. This
figure shows the computed clusters as well as the allocation function
for each element of the original distribution. Our second numerical
experiment is quantitative. In Table \ref{tab:2998}, we present standard
metric computations for a typical run of these algorithms with a high
number of clusters. This table illustrates the fact that each algorithm
is specialized in minimizing a particular functional, such as inertia
for k-means and the Maximum Mean Discrepancy (MMD) for kernel
clustering.

In view of these numerical experiments, the following comments are in
order.

\begin{itemize}
\item
  The first method is the greedy discrepancy algorithm \eqref{SDSD}.
  Observe that Figure \ref{fig:ccm} shows that the first clusters are of
  rather poor quality. However, Table \ref{tab:2998} exhibits an
  excellent ratio performance / execution time. Indeed, we regard this
  method as being statistical in nature, which usually performs better
  with a high number of clusters.
\item
  The second method is the sharp discrepancy algorithm \eqref{SDS}. This
  method usually gives the best results for our purpose; however Table
  \ref{tab:2998} shows that it comes at the expense of having a long
  running time.
\item
  The third method is the so-called scikit's k-means (mini batch)
  algorithm. It is a typical example where the k-means algorithm
  provides very good results, because we are studying a distribution
  whose modes are balanced and contain the same number of masses.
  Observe that cluster centers of k-means and sharp discrepancy differ
  significantly, resulting from different metric minimization.
\item
  The fourth method is the balanced clustering method \eqref{balanced},
  in which we randomly selected five points in \(X\) as our clusters
  \(Y\), and chose the discrepancy distance \(d_k(x,y)\). Here, we
  wanted to illustrate that this method is very stable, quite flexible,
  and Table \ref{tab:2998} also shows that is is numerically very
  efficient.
\end{itemize}

\begin{longtable}[t]{l|r|r|r|r|r}
\caption{\label{tab:2998}Supervised algorithm performance indicators}\\
\hline
method & $N_X$ & $N_Y$ & exec. time & inertia & MMD\\
\hline
\endfirsthead
\caption[]{Supervised algorithm performance indicators \textit{(continued)}}\\
\hline
method & $N_X$ & $N_Y$ & exec. time & inertia & MMD\\
\hline
\endhead
greedy & 1024 & 128 & 0.0389113 & 1028.1807 & 0.0000253\\
\hline
sharp disc. & 1024 & 128 & 0.5123663 & 1955.1476 & 0.0000102\\
\hline
kmeans & 1024 & 128 & 0.1977298 & 408.9355 & 0.0007832\\
\hline
bal. random & 1024 & 128 & 0.0190043 & 979.9428 & 0.0025103\\
\hline
\end{longtable}

\subsection{Supervised learning with
clustering}\label{supervised-learning-with-clustering}

\begin{figure}
\centering
\includegraphics{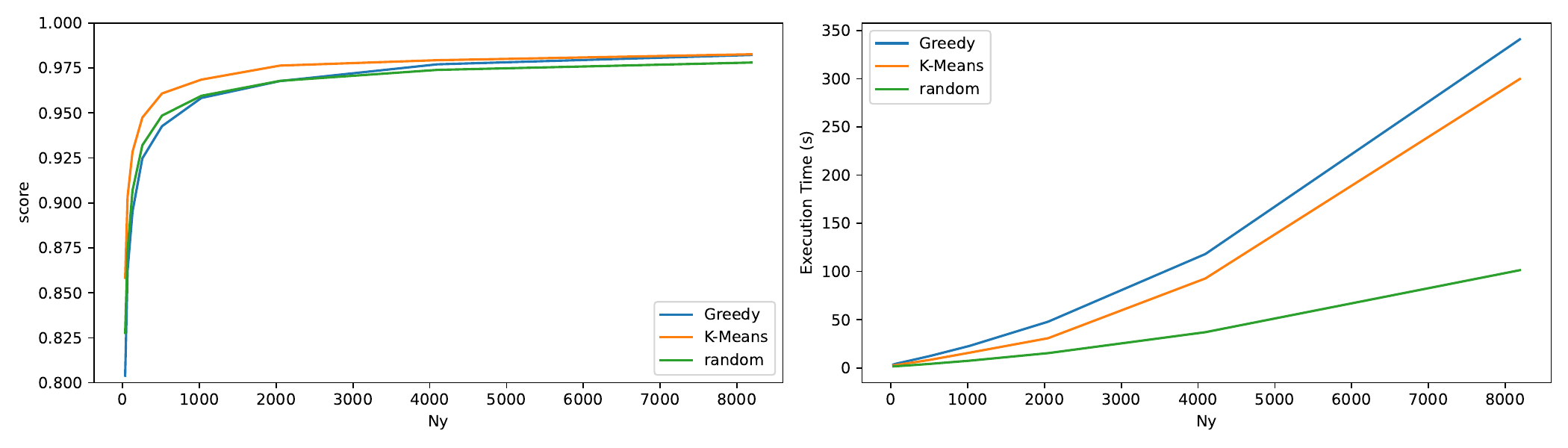}
\caption{Performance of interpolations for a given kernel - MNIST
problem}\label{fig:xy}
\end{figure}

\begin{figure}
\centering
\includegraphics{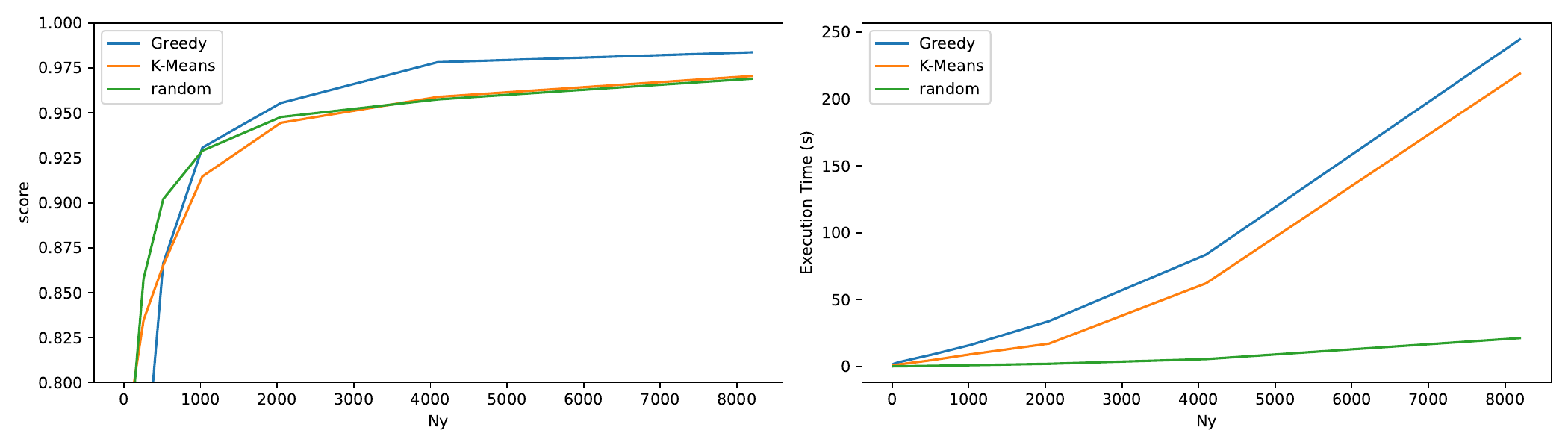}
\caption{Performance of extrapolations for a given kernel - MNIST
problem}\label{fig:yy}
\end{figure}

Next, we turn our attention to investigating the behavior of the
operator \eqref{FIT}, without regularization, for the interpolation mode
\(Y \subset X\) in Figure \ref{fig:xy}, as well as the extrapolation
mode \(Y = X\) in Figure \ref{fig:yy}, while \(Y\) is computed by some
clustering method.

In view of these results, we then replace our kernel by the Matern-based
kernel, which is known to be more robust in high dimensions, namely:
\(k(x,y) = \exp(|x-y|_1)\) with a scaling map \(k\circ S\), with
\(S=S_1\circ S_2\circ S_3\). Here, \(S_3(\cdot)\) maps the input data
\(X\) to the unit cube, \(S_2(\cdot)= \text{erf}^{-1}(\cdot)\)
(\(\text{erf}\) being the standard error function), and
\(S_1(\cdot) = \frac{\cdot}{\alpha}\) with scale factor
\(\alpha=\frac{1}{N_XN_Y}\sum_{n,m} |x^n-y^m|_1\).

Our numerical test considers a typical academic example of images
classification, called the MNIST problem (``Modified National Institute
of Standards and Technology''). This test is quite interesting in the
context of the present paper since the size of the training set is too
large in order to be treated directly with the operator \eqref{FIT} in
extrapolation mode, but is sufficiently small to allow to perform
various experiments. The MNIST dataset contains \(N_X=60,000\) training
images and \(N_Z=10,000\) testing images. The dimension of this problem
is \(D=784\), each image being coded with \(28 \times 28\) pixels. There
are \(10\) possible digits, from \(0\) to \(9\). Since this is a
classification problem, we used the probabilistic version of the
operator \eqref{FIT} given by \[
f_{k,\theta}(\cdot) = \text{softmax} (\log(\tilde{f})_{k,\theta})(\cdot),
\] where the softmax function is
\(\text{softmax}(\cdot) = \frac{\exp(\cdot)}{|\exp(\cdot)|_1}\) and
\(\tilde{f}\) is a projection of \(f\) into the unit cube obtained by
the proportional fitting algorithm. In Figure \ref{fig:xy} (left-hand),
we present the MNIST score obtained, that counts the number of correct
predictions of a predictive method as follows
\(\frac{\# \{ f_{k,\theta}(Z) = f(Z)\}}{N_X}\). Execution time are shown
in the right-hand plot. In these figures, \(N_y\) corresponds to the
number of computed clusters.

Our comments are as follows

\begin{itemize}
\tightlist
\item
  Comparing Figure \ref{fig:xy} (left-hand) shows that clustering
  methods achieves better scores than random picking, however they
  deteriorate the ratio score versus execution time. With this
  interpolation approaches, all our tests measured scores under 98 \%
  accuracy.
\item
  The best scores and ratio score versus execution time are obtained in
  Figure \ref{fig:yy} by the greedy clustering method. This method can
  reach 98.4 \% accuracy with as few as 2000 well-selected samples among
  the 60000, and adding more samples do not improve the results. Indeed,
  this score is the maximum recorded by a Gaussian kernel for this test,
  and other, convolutional kernels should be used to reach higher
  scores. We say that this kernel is \textit{saturating} with training
  set of greater size, an effect that is reminiscent of over-fitting
  with other methods.
\item
  The method sharp discrepancy was disqualified for this test, as
  execution time of this method was prohibitive for such dataset sizes
  and cluster numbers.
\end{itemize}

\subsubsection{Multiscale supervised
learning}\label{multiscale-supervised-learning}

We now illustrate the behavior of the multi-scale method \eqref{MS} for
the MNIST problem, reporting as previously scores and execution times.
For this experiment we selected \(N_Y=[5,\ldots, 80]\) clusters with
various clustering methods, running the MNIST test in full extrapolation
mode on each clusters. As it is paramount to retrieve balanced clusters,
we ran the balanced clustering method \eqref{balanced} and computed the
assignation function for the computed clusters, and also included random
picking clusters to have a benchmark.

\begin{itemize}
\tightlist
\item
  Considering scores, no method seems to pop up from this experience.
  Indeed, the balancing method \eqref{balanced} erased all differences,
  and even random picking achieves respectable scores in this test.
\item
  The less number of clusters comes with better scores and higher
  execution time. With five clusters, we are close to kernel saturation.
\item
  Note that the sharp discrepancy algorithm necessitates to compute the
  full Gram matrix at initialization. This sole operation took more than
  100 seconds and is responsible for the computational cost of this
  algorithm.
\end{itemize}

\begin{figure}
\centering
\includegraphics{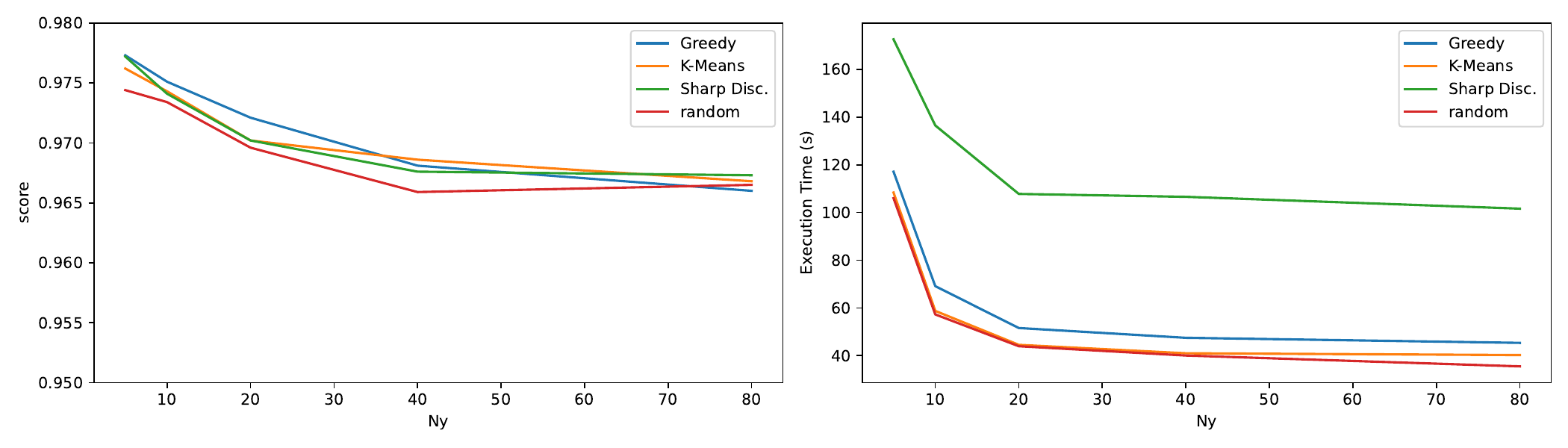}
\caption{Performance of multiscale learning for a given kernel - MNIST
problem}
\end{figure}

This experiment showed that this multiscale method indeed trades
accuracy versus time, but with sustainable losses, and is very
permissive with the cluster choices. This robust method allow to address
predictive kernel-based methods on large scale dataset, a single
computing unit as the CPU of a laptop being able to handle a bandwith
roughly between 1000 to 10000 samples by seconds, depending on the
problem.

\subsection{Optimal transport}\label{optimal-transport}

\subsubsection{Numerical comparison}\label{numerical-comparison}

As we suggest in this paper a different approach to optimal transport
than the mainstream one, that is the Sinkhorn algorithm with
regularization, we propose in this paragraph a benchmark of both
approaches, tailored to our needs. To that aim, we used the very same
methodology as the one in Pooladian and Niles-Weed \cite{NiPa:2021}
\footnote{we also used the code of these authors, kindly made publicly available.},
that is the following: consider a distribution \(X \in \mathbb{R}^d\),
compute \(Y = S(X)\) for a map having form \(S=\nabla h\), \(h\) convex,
so that \(S\) is an optimal transport map, scramble \(Y\) and input
\(X,Y\) into any method capable to determine \(S_{k,\theta}(\cdot)\).
The test set \(Z\) to benchmark with is defined sampling from the same
distribution as \(X\), and we measure the mean square error
\(|S_{k,\theta}(Z) - S(Z)|_{\ell^2}\).

\begin{figure}
\centering
\includegraphics{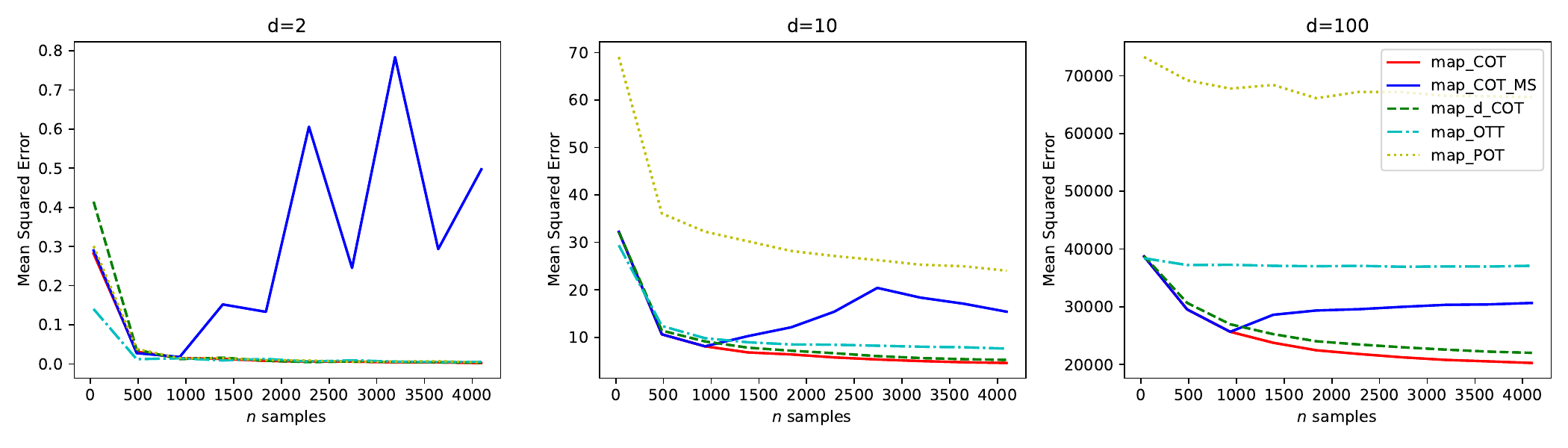}
\caption{Performance of LSAP vs Sinkhorn Optimal
Transport}\label{fig:otperf}
\end{figure}

\begin{figure}
\centering
\includegraphics{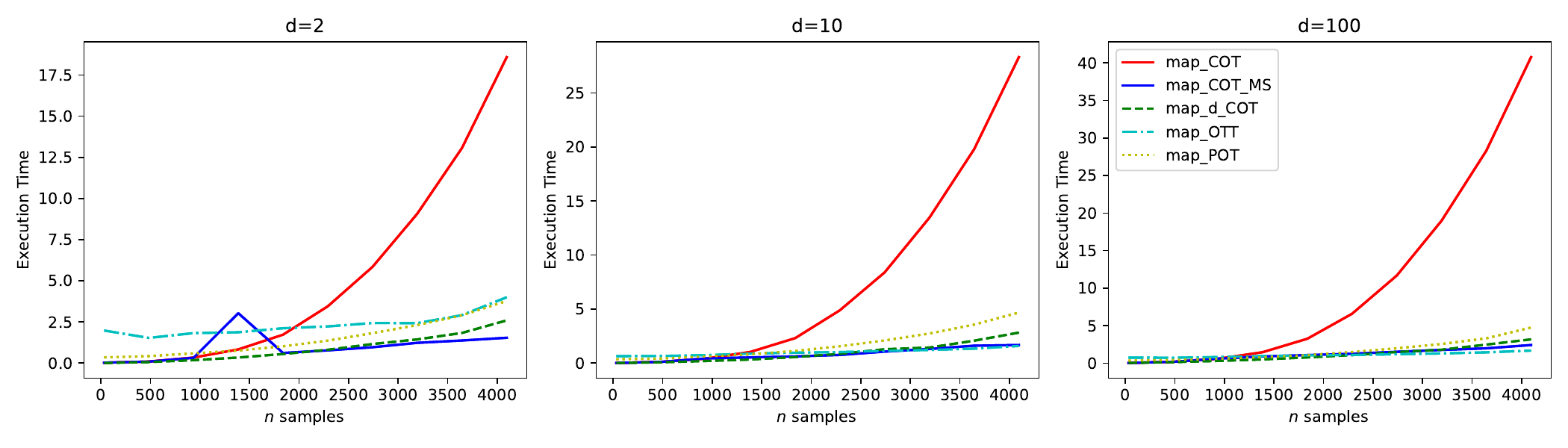}
\caption{Execution time of LSAP vs Sinkhorn Optimal
Transport}\label{fig:ottime}
\end{figure}

We choose \(X\) as the uniform distribution, \(S(X)=X|X|_2^2\), and
reported in Figure \ref{fig:otperf} the mean square error for different
dimension \(d=[2,10,100]\), different training set size
\(N=[256,512,...,4096]\), as well as corresponding execution time in
Figure \ref{fig:ottime} for different methods. We choose to benchmark
five methods, as follows.

\begin{itemize}
\tightlist
\item
  Three methods COT, COT Parallel, COT MS are based upon \eqref{ED}. The
  first COT uses an exact LSAP solver, and is not parallel. COT Parallel
  is a parallelized, sub-optimal solver for LSAP. COT MS is a multiscale
  method using \eqref{KM}, set arbitrarily with \(C = N/256\) clusters,
  in order to check that this method produces the same result as COT for
  the first example \(N=256\).
\item
  Two methods, POT and OTT, are based on a Sinkhorn regularization of
  optimal transport, using the publicly availables OTT and POT
  libraries. The Sinkhorn method is quite sensitive to regularization
  parameters, so we had to guess the best parameters, that is the
  smallest one for which the Sinkhorn algorithm converges in our case.
  OTT provides a method to guess this parameter, and we use trial and
  error to tune this paprameter for POT.
\end{itemize}

Our conclusions are as follows.

\begin{itemize}
\item
  Performance: A pure combinatorial method, being an exact approach,
  outperforms regularization-based methods such as Sinkhorn in all
  cases, as expected.
\item
  Scalability: The trade-off is scalability. As shown in Figure
  \ref{fig:ottime}, the POT method is faster and can run in parallel on
  GPUs, making it highly efficient for larger datasets.
\item
  Asymptotic Efficiency: The fastest asymptotic method is COT MS, which
  is linear in the dataset size and computes exact transport. However,
  this method also trades accuracy for computational time, particularly
  in lower-dimensional cases.
\end{itemize}

Indeed, the purpose of Sinkhorn algorithms ---------and regularized
optimal transport approaches in general--------- is to address
regularized optimal transport problems, which are desirable in certain
applications. These algorithms are efficient and scalable for such tasks
but cannot compete with direct methods for non-regularized optimal
transport problems, particularly in the context of small to medium-sized
datasets, as considered in this numerical experiment.

\subsubsection{Remark on the multiscale optimal transport
approach}\label{remark-on-the-multiscale-optimal-transport-approach}

In the previous experiment, we benchmarked a multiscale approach for
Optimal Transport between a source distribution \(X\) and a target
distribution \(Y\). For large distribution sizes \(N\), computing a
linear sum assignment (LSAP) problem directly between \(X\) and the
target can be computationally prohibitive. To address this, the approach
involves splitting the source and target distributions into \(C\)
clusters using balanced clustering. Each cluster is treated as a
manageable batch, enabling OT computations to be performed independently
and in parallel for each batch. The results from these batches are then
combined to reconstruct the estimated target distribution. This is a
fast method, but it relies on an approximation of the transport map that
we illustrate qualitatively in this paragraph.

\begin{figure}[H]

{\centering \includegraphics[width=0.6\linewidth,]{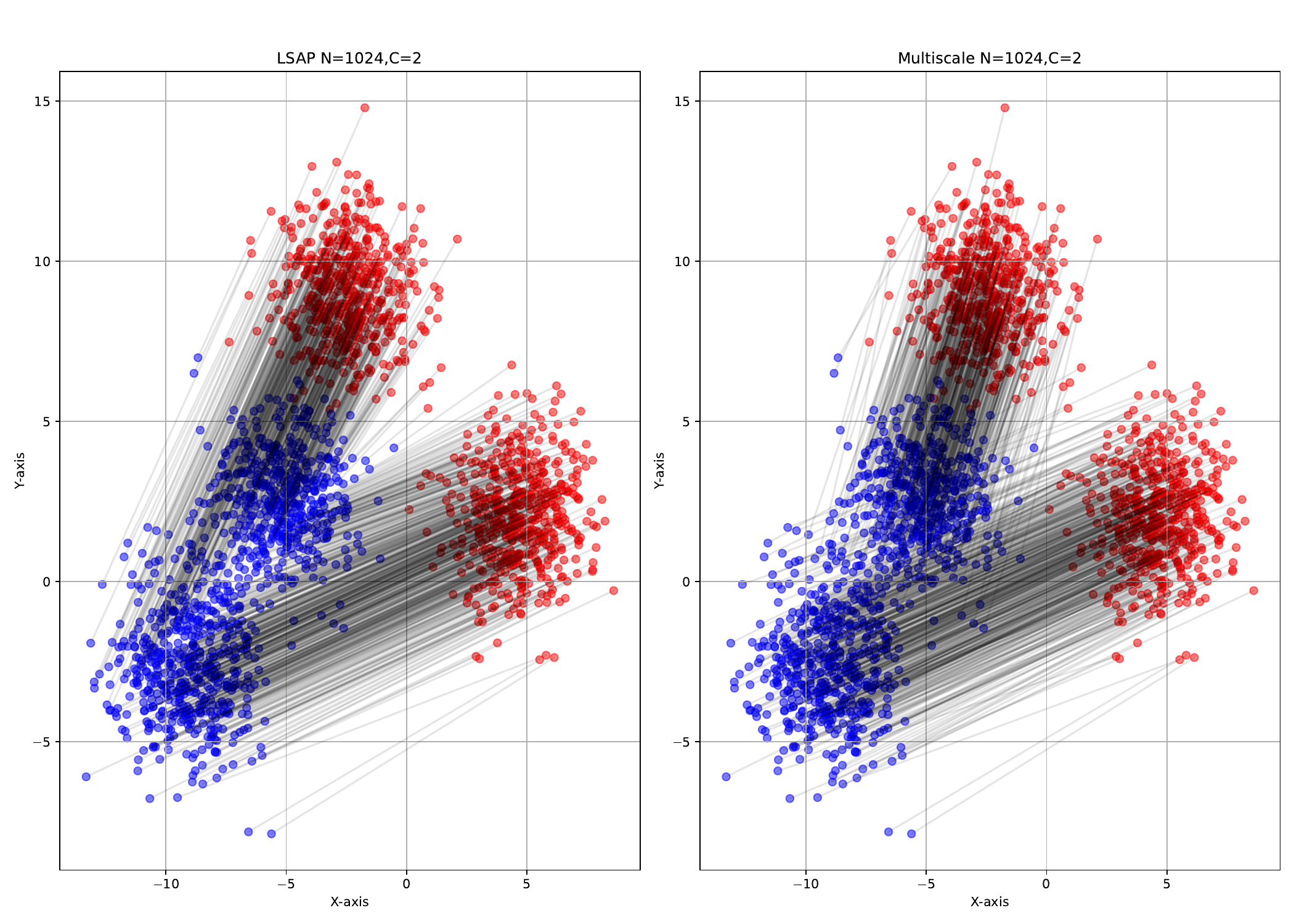} 

}

\caption{\label{fig:LSAPOT}Difference of assignments between LSAP and mutiscale OT.}\label{fig:unnamed-chunk-4}
\end{figure}

In this experiment, the source \(X\) and the target \(Y\) are sampled
from Gaussian mixtures with specified means and variances. The
multiscale optimal transport framework maps \(X\) to \(Y\) using a given
number of clusters \(C = 2\). We consider a small number \(N=1024\) and
the LSAP as a reference value, matching \(X\) to \(Y\) directly based on
the pairwise euclidean distances. A multiscale optimal transport is then
trained on the resulting assignment, based on the same distance. For
both approaches, we plot in figure \ref{fig:LSAPOT} the distributions,
\(X\) in red, \(Y\) in blue, and each black line gives the assignment.
The comparaison between the two figure allows to appreciate the
approximation made by this approach: when using the Euclidean distance,
assignments lines should never cross eachothers, while the multiscale
method do, as illustrated by these two figures.

\subsection{Conditional sampling}\label{conditional-sampling}

We refer to \cite{LMM-Wilmott} for some examples of application to
time-series forecasting with the conditional sampling algorithm
\eqref{condalg}. In this paragraph, we chose to illustrate our
conditional image generator in modifying specific attributes in images
from the CelebA dataset. This dataset is a collection of images having a
RGB resolution of \(178\times 218 \times 3=116412\). Each image comes
with 40 binary labels indicating facial attributes like hair color,
gender and age.

\begin{figure}[H]
\includegraphics[width=1\linewidth,]{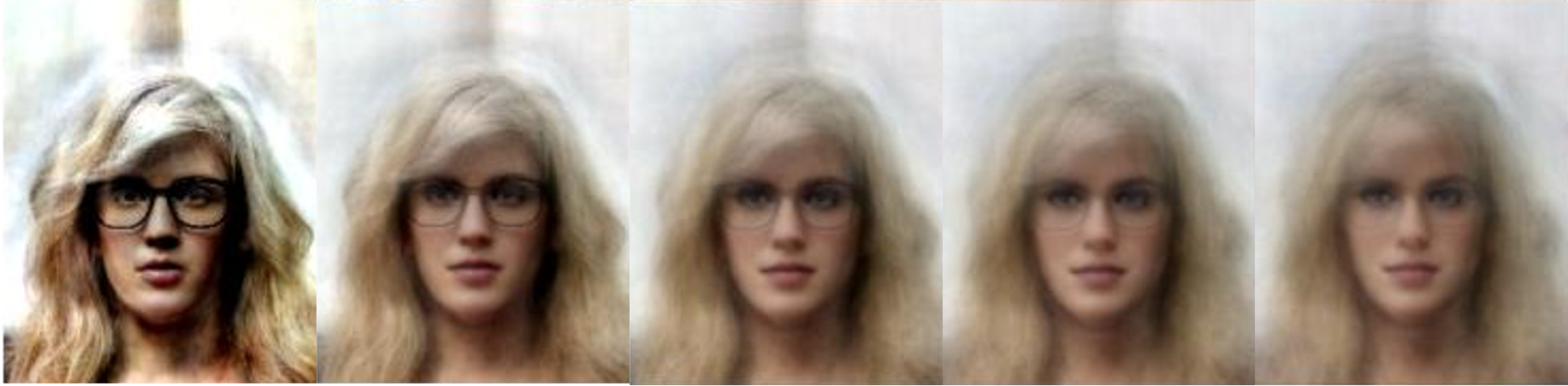} \caption{\label{fig:celeba} Progressive removal of the 'glasses' attribute from a representative image in the CelebA dataset. The first column shows the original image with glasses, and each subsequent column represents an incremental reduction in the 'glasses' attribute, eventually reaching the target state of 'no glasses' in the final column. This illustrates the generator's capacity to modify specific attributes while maintaining overall facial features, with slight blurring effects appearing as the glasses are gradually removed.}\label{fig:unnamed-chunk-5}
\end{figure}

Considering these datas, we use the algorithm \eqref{condalg} to model
the conditional distribution \(Y|X\) with \(Y \in \mathbb{R}^{116412}\)
and \(X \in \mathbb{R}^{40}\), considering a subset of 1,000 images,
randomly selected from those with the attributes {[}Woman, light
makeup{]}. From this subset, we identified images with the attribute
{[}hat,glasses{]} (denoted as {[}+1,+1{]}) and select a single
representative image for illustration, shown in the first column of
Figure \ref{fig:celeba}. Our objective is to progressively remove the
hat and glasses attribute, reaching the target state {[}-1,-1{]} (i.e.,
no hat, no glasses).

The latent variable, denoted \(\eta_y\) in \eqref{condalg}, is modeled
as a standard Gaussian distribution in 25 dimensions, and we do not use
a latent variable for the attributes, i.e.~\(\eta_X=X\). While all other
latent components of the selected image remain constant,
{[}hat,glasses{]} attribute are gradually transitioned from {[}+1{]} to
{[}-1{]} in increments of 0.4 across each column in Figure
\ref{fig:celeba}. The final column represents the target state, where
the subject appears without glasses and hat. This experiment assesses
the generator's ability to exactly reproduce the training set, and to
smoothly interpolate attribute changes to produce images that coherently
reflect the conditioned modifications.

The dimensionality of the latent space significantly influences the
results. A low-dimensional latent space tends to retain the appearance
of glasses, while a high-dimensional space successfully removes this
feature but may result in overly blurred and less distinct images. We
manually tuned this parameter to achieve a balance between attribute
removal and image clarity in the generated figure.

\subsection{The Bachelier problem}\label{the-bachelier-problem}

We present the so-called Bachelier test, a pertinent test for
mathematical finance applications, in order to test transition
probability algorithms, already used in \cite{HG:2020} to provide a test
for a neural network approach to transition probability algorithms, and
also in \cite{LMM-SSRN1} to benchmark a variant of Algorithm
\ref{sumalg}.

We benchmark \eqref{sumalg} to others algorithms, namely the Sinkhorn
method, and the Nadaraya Watson method, providing an alternative method
to compute conditional expectations. The test is described as follows.

\begin{itemize}
\item Consider a Brownian motion $t \mapsto X_t \in \mathbb{R}^D$, satisfying $dX_t = \sigma dW_t$,
where the matrix $\sigma \in \mathbb{R}^{D,D}$ is randomly generated. The initial condition is $X_0 = 0$ w.l.o.g. Let $\omega \in \mathbb{R}^D$, randomly generated, satisfying $|\omega|_1 = 1$ and denote the basket values $b_t = <\omega,X_t>$. This last process follows a univariate Brownian motion $d b_t = \theta dW_t$. We normalize $\sigma$ in order to retrieve a constant value for $\theta$, fixed to 0.2 in our tests.

\item Consider two times $1 = t_1 < t_2 = 2$, $t_2$ being the maturity of an option, that is a function
denoted $P(x) = \max(b(x) - K, 0)$. The goal of this test is to benchmark some methods aiming to compute the conditional expectation $\mathbb{E}^{X_{t_2}}[P(\cdot) \mid X_{t_1}]$, for which the reference value is given by the so-called Bachelier formula
\begin{equation}\label{BF}
  f(\cdot):=\mathbb{E}^{X_{t_2}}[P(\cdot) \mid X_{t_1}] = \theta \sqrt{t_2-t_1} \ p(d) + (b_{t_1}-K)c(d),\quad d=\frac{b_{t_1}-K}{\theta\sqrt{t_2-t_1}},
\end{equation}
where $p$ (resp. $c$) holds for the cumulative (resp. density) of the normal law.

\item Considering two integer parameters $N,D$, in this test we consider three samples of the Brownian motion $X \sim X_{t_1},Y \sim X_{t_2}, Z\sim X_{t_1}$ in $\mathbb{R}^{N,D}$, and produce a transition probability to approximate the transition probability matrix $\Pi(Y | X)$. We finally output the mean square error $\text{err}(Z) = \|f(Z)- \Pi(Y | Z) P(Y) \|_{\ell^2}$.
\end{itemize}

\begin{figure}
\centering
\includegraphics{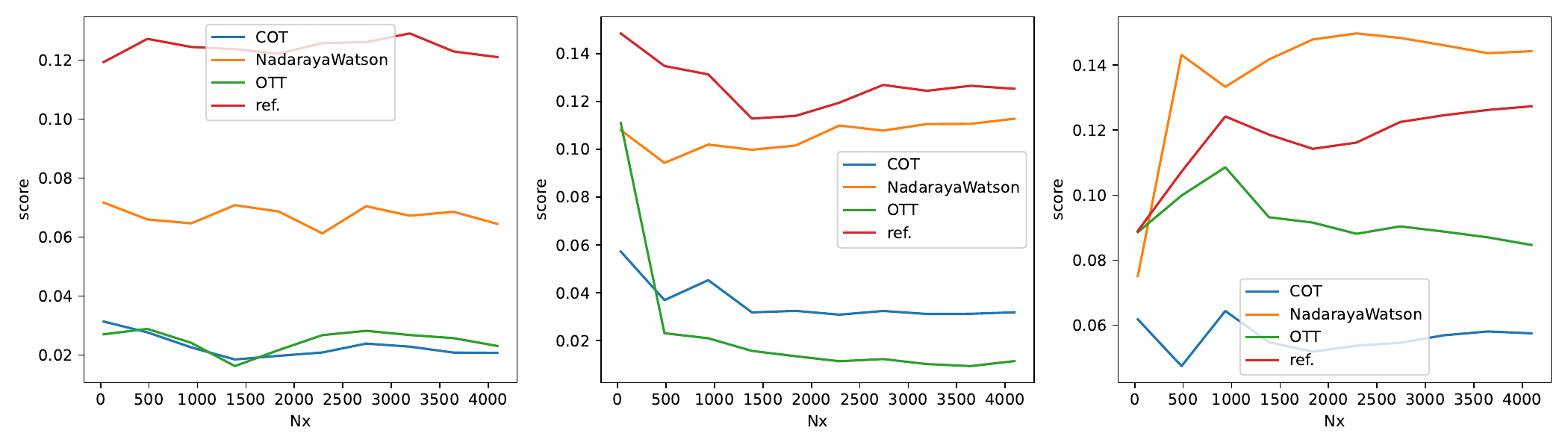}
\caption{Benchmark of scores}\label{fig:BacPerfs}
\end{figure}

\begin{figure}
\centering
\includegraphics{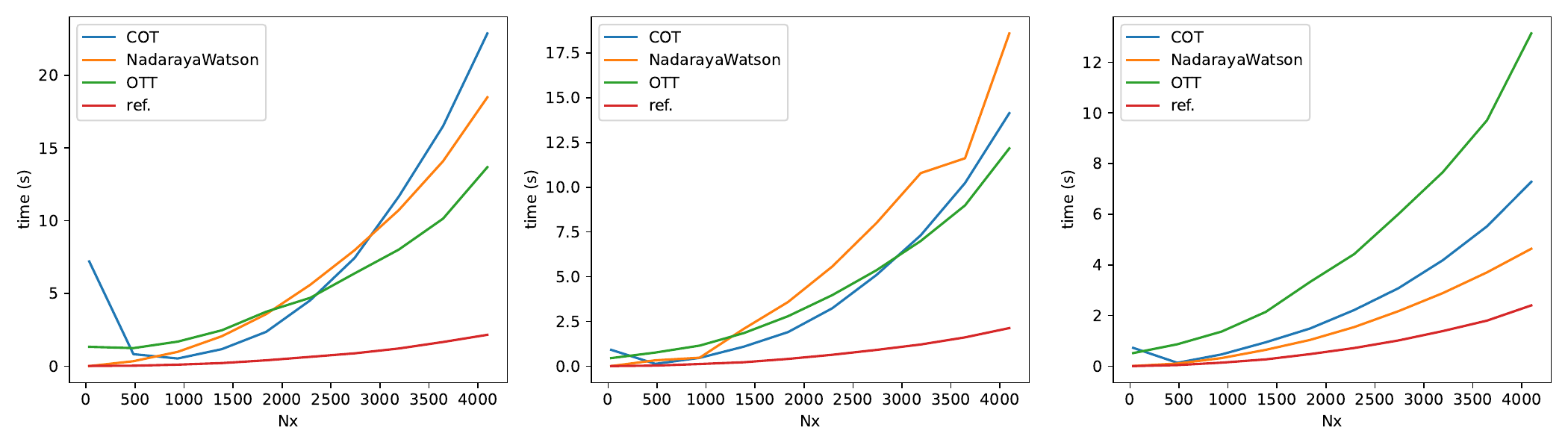}
\caption{Execution time}\label{fig:ottime-deux}
\end{figure}

Figure \ref{fig:BacPerfs} (resp. Figure\ref{fig:ottime-deux}) shows the
score
\(\frac{\text{err}(Z)}{\|f(Z)\|_{\ell^2}+\|\Pi(X_{t_2} | Z) P(X_{t_2})\|_{\ell^2}}\)
for various choice of \(N\) and \(D\) (resp. the execution time in
seconds) for four methods.

\begin{itemize}
\item COT illustrates our results with Algorithm \ref{sumalg}.
\item OTT consists in our best trial using Sinkhorn algorithm with the OTT library.
\item Ref is a reference, naive value, computing $P_{k,\theta}(Z)$ with the interpolation formula \eqref{FIT}.
\item Nadaraya Watson implements \eqref{NW}.
\end{itemize}

All these methods (labelled with an index \(m\)) output a transition
probability matrix \(\Pi^m(Y|X)\), which implies the estimation
\(f^m(X) = \Pi(Y|X) P(Y)\). We used the extrapolation method \eqref{FIT}
in order to extrapolate \(f^m_{k,\theta}(Z)\) for all methods. In this
test, we observed that both the Sinkhorn algorithm and \ref{sumalg}
performed well across a wide range of parameters \(N\) and \(D\). Both
approaches reliably computed the transition probability matrix. Notably,
the naive extrapolation method exhibited surprisingly good performance
in the high-dimensional case (\(D=100\)). This unexpected result raises
questions about the relevance of this test in high dimensions, where a
simple linear regression method could achieve similar accuracy when
approximating the Bachelier formula \eqref{BF}.

\begin{thebibliography}{00} 

\bibitem{Bach}
{\sc J. Altschule and F. Bach,}
Massively scalable Sinkhorn distances via the Nystrom method, available at https://arxiv.org/pdf/1812.05189

\bibitem{BTA}
{\sc A. Berlinet and C. Thomas-Agnan,}
{\it Reproducing kernel Hilbert spaces in probability and statistics,}
Springer Science, Business Media, LLC, 2004. 

\bibitem{Brezis:2018} 
{\sc H. Brezis,}
Remarques sur le probl\`eme de Monge–Kantorovich dans le cas discret, 
Comptes Rendus Math. 356 (2018), 207--213.

\bibitem{Cuturi:2013} 
{\sc M. Cuturi,} 
Sinkhorn distances: lightspeed computation of optimal transport, 
 Advances in Neural Information
Processing Systems 26: 27th Annual Conference on Neural Information Processing Systems 2013, C.J.C. Burges, L. Bottou, Z. Ghahramani, and K.Q. Weinberger, editors, 
Proceedings of a meeting held December 5-8, 2013, Lake Tahoe, Nevada, United States, pp. 2292--2300.

\bibitem{Bachdeux} 
{\sc N. Doum\'eche, F. Bach, G. Biau, and C. Boyer,} 
Physics-informed kernel learning, available at https://arxiv.org/abs/2409.13786

\bibitem{EFO} 
{\sc F. Eckerli and J. Osterrieder,} 
Generative adversarial networks in finance: an overview, 2021, available as ArXiv:2106.06364.
 
\bibitem{García:2024} 
{\sc A. García},
Greedy algorithms: a review and open problems, 2024, available as ArXiv:2408.08935.

\bibitem{GBRS}
{\sc A. Gretton, K.M. Borgwardt, M. Rasch, B. Sch\"{o}lkopf, and A.J. Smola,}
A kernel method for the two sample problems, 
Proc. 19th Int. Conf. on Neural Information Processing Systems, 2006, pp.~513--520.

\bibitem{HT:1985} 
{\sc J.J.Hopfield and D.Tank,}
'Neural' computations of decisions in optimization problems, Biological Cybernetics 52 (1985), 141––152.

\bibitem{HG:2020} 
{\sc B.N. Huge and A. Savine,}
Differential machine learning, January 2020. Available at http://dx.doi.org/10.2139/ssrn.3591734

\bibitem{LM-CRAS}
{\sc P.G. LeFloch and J.-M. Mercier,}
A new method for solving Kolmogorov equations in
mathematical finance, C. R. Math. Acad. Sci. Paris 355 (2017), 680–686.

\bibitem{PLF-JMM-estimate}
{\sc P.G. LeFloch and J.-M. Mercier,} 
Mesh-free error integration in arbitrary dimensions: a numerical study of discrepancy functions,
Comput. Methods Appl. Mech. Engrg. 369 (2020), 113245.

\bibitem{PLF-JMM-Wilmott}
{\sc P.G. LeFloch and J.-M. Mercier,} 
The transport‐based mesh‐free method: a short review,
 Wilmott Journal 109 (2020), 52––57.

\bibitem{PLF-JMM:2023}
{\sc P.G. LeFloch and J.-M. Mercier,} 
A class of mesh-free algorithms for some problems arising in finance and machine learning,
J Sci Comput 95 (2023), 75. https://doi.org/10.1007/s10915-023-02179-5

\bibitem{LMM-SSRN1} 
{\sc P.G. LeFloch, J.-M. Mercier, and S. Miryusupov,}
A kernel based method for computing conditional expectations. (March 29, 2021). Available at SSRN: https://ssrn.com/abstract=3814704 or http://dx.doi.org/10.2139/ssrn.3814704

\bibitem{LMM-Wilmott} 
{\sc P.G. LeFloch, J.-M. Mercier, and S. Miryusupov,}
Extrapolation and generative algorithms for three applications in finance, 
The Willmot Journal, September 2024, pp. 54--60. Available at  https://arxiv.org/abs/2404.13355

\bibitem{LMM} 
{\sc P.G. LeFloch, J.-M. Mercier, and S. Miryusupov,}
CodPy: A Python Library for Machine Learning, Mathematical Finance, and Statistics, monograph in preparation available at http://arxiv.org/abs/2402.07084.
Code available at \url{https://pypi.org/project/codpy/}

\bibitem{LMMprep} 
{\sc P.G. LeFloch, J.-M. Mercier, and S. Miryusupov,}
Physics-informed learning with reproducing kernels, in preparation.

\bibitem{MaFr:2014}
{\sc M.I. Malinen and P. Fränti},
Balanced K-means for clustering. Structural, syntactic, and statistical pattern recognition, Lecture notes in computer science, Vol. 8621, pp. 32--41.

\bibitem{Falkon:2020}
{\sc G. Meanti, L. Carratino, L. Rosasco, Lorenzo, and A. Rudi,}
Kernel methods through the roof: handling billions of points efficiently (2020), 
Advances in Neural Information Processing Systems 32. Available as https://arxiv.org/abs/2006.10350

\bibitem{MF2011} 
{\sc F. Memoli,} 
Gromov–Wasserstein distances and the metric approach to object matching,
Foundations of Comp. Math. 11 (2011), 417--487.

\bibitem{Cuturi:2016} 
{\sc G. Peyré, M. Cuturi, and J. Solomon,} 
Gromov-Wasserstein averaging of kernel and distance matrices, 
Proc. Machine Learning Research 48 (2016), 2664––2672. 
 
\bibitem{NiPa:2021}
{\sc  A-A. Pooladian and J. Niles-Weed,}
Entropic estimation of optimal transport maps, available as https://arxiv.org/abs/2109.12004
 
\bibitem{B2001} 
{\sc B. Sch\"{o}lkopf, R. Herbrich, and A.J. Smola,}
A generalized representer theorem, 
in Computational learning theory, Springer Verlag, 2001,  pp. 416--426. 

\bibitem{SFVTP:2023} 
{\sc T. Séjourné, J. Feydy, F.-X. Vialard, A. Trouvé, and G. Peyré,}
Sinkhorn divergences for unbalanced optimal transport, 
available as https://arxiv.org/abs/1910.12958

\bibitem{TH:2019} 
{\sc E.D. Taillard and K. Helsgaun},
Popmusic for the travelling salesman problem, 
European Journal of Operational Research 272 (2019), 420--429.

\bibitem{Temlyakov:2011}
{\sc V. N. Temlyakov,}
{\it Greedy approximation,}
Vol 20, Cambridge University Press, 2011.

\bibitem{LiVlVe:2009}
{\sc G. Tzortzis and A. Likas,}
The global kernel k-means clustering algorithm,
2008 IEEE International Joint Conference on Neural Networks (IEEE World Congress on Computational Intelligence), Hong Kong, China, 2008, pp. 1977--1984.  

\bibitem{Villani:2009}
{\sc C. Villani,}
{\it Optimal transport, old and new,}
Springer Verlag, 2009. 

\bibitem{BLWY2}
{\sc K. A. Wang, X. Bian, P. Liu and D. Yan,}
DC2: a divide-and-conquer algorithm for large-scale kernel learning with application to clustering, 
2019 IEEE International Conference on Big Data, Los Angeles, CA, USA, 2019, pp. 5603--5610. Available at https://ieeexplore.ieee.org/document/9006565 

\bibitem{WaSaHa:2023} 
{\sc T. Wenzel, G. Santin, and B. Haasdonk},
Analysis of target data-dependent greedy kernel algorithms: convergence rates for f, fP, and f/P greedy, Constructive Approximation 57 (2023), 45--74.

\end{thebibliography}
\end{document}